\documentclass[12 pt]{amsart}

\usepackage{graphicx}
\usepackage{color}
\usepackage{amssymb, amsmath, amsthm}
\usepackage{mathptmx}

\usepackage{epsfig}
\usepackage{epstopdf}
\usepackage{graphicx}

\theoremstyle{plain}

\newtheorem{theorem}{Theorem}
\newtheorem{lemma}{Lemma}

\newtheorem{question}{Question}

\newcommand{\R}{\mathbb R} %REALS

\newcommand{\Q}{\mathbb Q} %RATIONALS
\newcommand{\Z}{\mathbb{Z}}

\newcommand{\bi}{\begin{itemize}}
\newcommand{\ei}{\end{itemize}}
\newcommand{\be}{\begin{enumerate}}
\newcommand{\ee}{\end{enumerate}}

\newcommand{\n}{\beta}

\newcommand{\emp}{\emptyset}
\newcommand{\X}{\times}

\newcommand{\eps}{\epsilon}

\newcommand{\A}{\alpha}
\newcommand{\pd}{\partial}

\newcommand{\Sau}{\mathcal{S}}
\newcommand{\FF}{\mathcal{F}}

\newcommand{\ceta}{\overline{\eta}}

\numberwithin{definition}{section}
\numberwithin{example}{section}
\numberwithin{lemma}{section}
\numberwithin{theorem}{section}
\numberwithin{corollary}{section}

\begin{document}
\title{Bridge and pants complexities of knots}
\author{Alexander Zupan}

\maketitle

\begin{abstract}
We modify an approach of Johnson \cite{johnson} to define the distance of a bridge splitting of a knot $K$ in a 3-manifold $M$ using the dual curve complex and pants complex of the bridge surface.  This distance can be used to determine a complexity, which becomes constant after a sufficient number of stabilizations and perturbations, yielding an invariant of $(M,K)$.  We also give evidence toward the relationship between the pants distance of a bridge splitting and the hyperbolic volume of the exterior of $K$.
\end{abstract}

\section{Introduction}
In the past decade, there has been great interest in studying topological properties of knots and manifolds via certain cell complexes associated to splitting surfaces.  In particular, the curve complex, pants complex, and dual curve complex have been employed in this regard.  The first such instance occurs in a paper of John Hempel \cite{hempel}, in which he uses the curve complex $C(\Sigma)$ of a Heegaard surface $\Sigma$ to define an integer complexity, the \emph{distance} (or \emph{Hempel distance}) $d(\Sigma)$ of $\Sigma$.  This distance refines the idea of strong irreducibility of Heegaard splittings and carries information about both essential surfaces contained in a 3-manifolds $M$ and alternate Heegaard splittings of $M$.  Haken's lemma implies that if $M$ is reducible, then for any Heegaard surface $\Sigma$ for $M$, $d(\Sigma) = 0$.  In his seminal paper, Hempel shows that if $M$ contains an essential torus, $d(\Sigma) \leq 2$ for all $\Sigma$, and Hartshorn generalizes this phenomenon, showing in \cite{hartshorn} that if $M$ contains an essential surface of genus $g$, then $d(\Sigma) \leq 2g$ for all splitting surfaces $\Sigma$.  In addition, Scharlemann and Tomova demonstrate that if $\Sigma'$ another Heegaard surface of genus $g'$ which is not a stabilization of $\Sigma$, then $d(\Sigma) \leq 2g'$ \cite{schartom}. \\

In the context of bridge splittings of knots, Bachman and Schleimer have used the arc and curve complex to adapt Hempel's distance to bridge splittings of knots, proving a result similar to that of Hartshorn: the distance of any splitting surface is bounded above by a function of the $\chi(S)$, where $S$ is an essential surface in the knot exterior $E(K)$ \cite{bacsch}.  Further, Tomova has proved that a distance similar to that of Bachman and Schleimer gives a lower bound on the genus of alternate bridge surfaces \cite{tomova}. \\

In \cite{johnson}, Johnson invokes the pants complex $P(\Sigma)$ and defines the related dual curve complex $C^*(\Sigma)$ of a Heegaard surface $\Sigma$ for a 3-manifold.  He proves that these complexes can be used to assign an integer complexity to $\Sigma$ and that this complexity converges to an integer $A(M)$ or $A^P(M)$ (depending on the complex used) upon taking a sequence of stabilizations of $\Sigma$.  The famous Reidemeister-Singer Theorem states that any two Heegaard surfaces have a common stabilization, and so $A(M)$ and $A^P(M)$ are invariants of $M$. \\

In this paper, we adapt Johnson's approach to define integer complexities $B(\Sigma)$ and $B^P(\Sigma)$ of a bridge splitting surface $\Sigma$ of a knot $K$ in a 3-manifold $M$.  We prove the following theorem:

\begin{theorem}
Let $K$ be a knot in a closed, orientable 3-manifold $M$, where $M$ has no $S^2 \X S^1$ summands, let $\Sigma$ be a splitting surface for $(M,K)$, and let $\Sigma^h_c$ be an $(h,c)$-stabilization of $\Sigma$.  Then the limits
\[ \lim_{h,c \rightarrow \infty} B(\Sigma) \qquad \text{ and } \qquad \lim_{h,c \rightarrow \infty} B^P(\Sigma)\]
exist.  Moreover, they do not depend on $\Sigma$ and thus define invariants $B(M,K)$ and $B^P(M,K)$ of the pair $(M,K)$.
\end{theorem}

The pants complex $P(\Sigma)$ of a surface $\Sigma$ is itself an interesting object of study; for instance, Brock has shown that it is quasi-isometric to the Teichm\"{u}ller space $T(\Sigma)$ equipped with the Weil-Petersson metric \cite{brock}, and Souto has revealed other surprising connections between pants distance and the geometry of certain 3-manifolds \cite{souto}.  We exploit a theorem of Lackenby \cite{lackenby} to prove the following:
\begin{theorem}\label{bounds}
Suppose $K$ is a hyperbolic 2-bridge knot, $\Sigma$ is a $(0,2)$-splitting surface for $K$, and $v_3$ is the volume of a regular hyperbolic ideal 3-simplex.  Then
\[ v_3(D_P(\Sigma) - 3) \leq \text{vol}(E(K)) < 10v_3(2D_P(\Sigma) - 3).\]
\end{theorem}

The paper proceeds in the following manner:  In Section \ref{stabpert}, we provide background information and include a proof of the analogue of the Reidemeister-Singer Theorem for bridge splittings in arbitrary manifolds.  In Section \ref{pnts} we introduce the curve, dual curve, and pants complexes, and in Section \ref{dist} we define the distance of a bridge splitting and prove several basic facts about this distance.  In Section \ref{BC}, we use the distance of the previous section to define bridge and pants complexity, and we prove the main theorem.  In Section \ref{prop}, we demonstrate several properties of these new invariants, and in Section \ref{crit}, we define the concept of a critical splitting in order to provide explicit calculations in Section \ref{calc}.  Finally, in Section \ref{hyp}, we discuss connections between pants distance and hyperbolic volume, and in Section \ref{quest} we include several interesting open questions.

\section{Stabilization and Perturbation}\label{stabpert}
We begin with definitions of compression bodies and Heegaard splittings.  Let $S$ be a closed surface.  A \emph{compression body} is defined as the union of $S \X I$ with a collection of 2-handles and 3-handles attached along $S \X \{0\}$, where the 3-handles cap off any 2-sphere boundary components.  We call $S \X \{1\}$ the \emph{positive boundary} of $C$, denoted $\pd_+ C$, and $\pd_-C = \pd C - \pd_+C$ the \emph{negative boundary} of $C$.  If $\pd_-C = \emp$, $C$ is a handlebody, and if $C = S \X I$, we call $C$ a \emph{trivial compression body}.  Let $M$ be a compact, orientable 3-manifold.  A \emph{Heegaard splitting} of $M$ is a decomposition of $M$ into two compression bodies $V$ and $W$, such that $V \cap W = \pd_+ V = \pd_+ W$.  We call $\Sigma = \pd_+V$ the \emph{Heegaard surface} of the splitting, and the genus of the splitting is defined to be the genus of $\Sigma$. \\

Analogously, a useful way to study a 1-manifold embedded in a 3-mani- fold is to decompose the space via a bridge splitting (from \cite{morsak}), for which we require the following definitions:  A collection of properly embedded arcs $\A$ contained in a compression body $C$ is called \emph{trivial} if each arc $\gamma \in \A$ either cobounds a disk with an arc contained in $\pd_+ C$ (such a disk is called a \emph{bridge disk}) or is properly isotopic to some arc $\{x\} \X I \subset S \X I$.  The first class of arcs is called \emph{$\pd_+$-parallel}; the second class, \emph{vertical}.  Let $K$ be a 1-manifold properly embedded in a compact, orientable 3-manifold $M$.  A \emph{bridge splitting} of $K$ is a decomposition of $(M,K)$ into the union of $(V,\A)$ and $(W,\n)$, such that $M = V \cup W$ is a Heegaard splitting, and $\A$ and $\n$ are collections of trivial arcs in $V$ and $W$ with $\pd \A \cap \pd_+ V = \pd \n \cap \pd_+ W$.  The surface $\Sigma = \pd_+ V$ is called a \emph{bridge surface}. \\

\begin{figure}[h]
  \centering
    \includegraphics[width=1\textwidth]{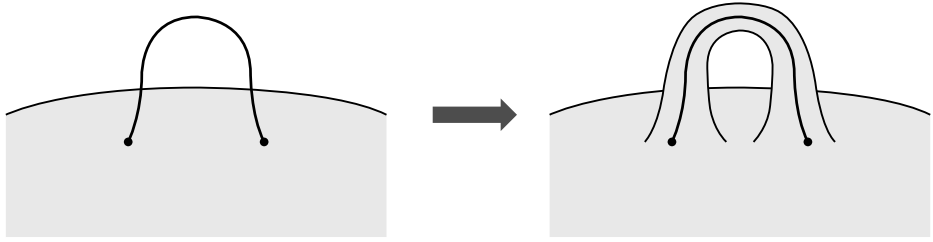}
    \caption{An elementary stabilization, in which $\eta(\gamma)$ is removed from $V$ and attached to $W$} \label{sta}
\end{figure}

\begin{figure}[h]
  \centering
    \includegraphics[width=1\textwidth]{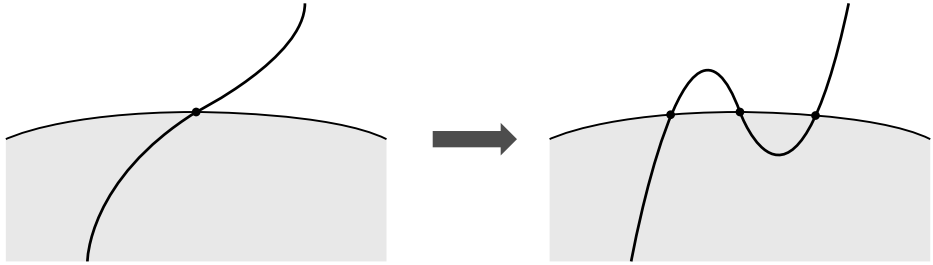}
    \caption{An elementary perturbation} \label{per}
\end{figure}

In general, we often refer to a bridge splitting by specifying only the bridge surface $\Sigma$, as $\Sigma$ uniquely determines $(V,\A)$ and $(W,\n)$.  Two bridge splittings $\Sigma$ and $\Sigma'$ are \emph{equivalent} if there is an isotopy of $(M,K)$ taking $\Sigma$ to $\Sigma'$.  Given a bridge splitting $(M,K) = (V,\A) \cup_{\Sigma} (W,\n)$, we can always make the splitting more generic in two ways.  To increase the genus of the splitting, let $\gamma \subset V$ be a $\pd_+$-parallel arc such that $\gamma \cap \A = \emp$, and let $\eta(\gamma)$ and $\overline{\eta}(\gamma)$ denote open and closed regular neighborhoods of $\gamma$ in $V$, respectively.  Define $V' = V - \eta(\gamma)$, $W' = W \cup \ceta(\gamma)$, and $\Sigma' = \pd_+ V' = \pd_+ W'$, so that $(M,K) = (V',\A) \cup_{\Sigma'} (W',\n)$ is a splitting of $(M,K)$ of higher genus.  This process is called \emph{elementary stabilization}.  On the other hand, we can perturb $K$ near some point of $K \cap \Sigma$ in order to add an extra trivial arc to $\A$ and $\n$.  This process is called \emph{elementary perturbation}.  See Figures \ref{sta} and \ref{per}.  If $K$ is a knot, then $\A$ and $\n$ consist only of $\pd_+$-parellel arcs, and so $|\A| = |\n|$.  In this case, we say that $(M,K) = (V,\A) \cup_{\Sigma} (W,\n)$ is a $(g,b)$-bridge splitting and $\Sigma$ is a $(g,b)$-bridge surface, where $g$ is the genus of $\Sigma$ and $b = |\A|$.  Note that elementary stabilization and elementary perturbation change a $(g,b)$-splitting into $(g+1,b)$- and $(g,b+1)$-splittings, respectively. \\

If $\Sigma^*$ is equivalent to an $(h,c)$-bridge surface obtained by applying some number of elementary stabilizations and perturbations to $\Sigma$, we say that $\Sigma^*$ is an $(h,c)$-\emph{stabilization} of $\Sigma$, and if $\Sigma$ and $\Sigma'$ have a common $(h,c)$-stabilization $\Sigma^*$, we say that $\Sigma$ and $\Sigma'$ are \emph{stably equivalent}.  Let $K$ be a system of trivial arcs contained in a handlebody $H$, so that every arc of $K$ is $\pd_+$-parallel.  The trivial bridge splitting of $(H,K)$ is $(V,\A) \cup_{\Sigma} (W,\n)$, where $W = \pd H \X I$ and contains only vertical arcs (forcing $(V,\A)$ to be homeomorphic to $(H,K)$).  We will employ the following special case of a theorem proved in \cite{hayashi} by Hayashi and Shimokawa:
\begin{theorem}\label{stand}
Suppose that $K$ is a system of trivial arcs contained in a handlebody $H$, with $(V,\A) \cup_{\Sigma} (W,\n)$ a bridge splitting of $(H,K)$ such that $V$ is a handlebody and $W$ is a trivial compression body.  Then $\Sigma$ is the result of some number (possibly zero) of elementary perturbations of the trivial splitting.
\end{theorem}
We use this fact in the next theorem, analogous to the Reidemeister-Singer Theorem which states that any two Heegaard splittings of a manifold $M$ are stably equivalent.  In fact, our proof is modeled on the recent proof of this famous theorem by Lei \cite{lei}.
\begin{theorem}\label{stab}
Let $K$ be a link in a closed 3-manifold $M$, with bridge splittings $(V,\A) \cup_{\Sigma} (W,\n)$ and $(V',\A') \cup_{\Sigma'} (W',\n')$.  Then $\Sigma$ and $\Sigma'$ are stably equivalent.  Moreover, a common stabilization $\Sigma^*$ can be obtained by applying a sequence of elementary stabilizations followed by a sequence of elementary perturbations to $\Sigma$ and $\Sigma'$.
\begin{proof}
First, the Reidemeister-Singer Theorem implies that the underlying Heegaard splittings $V \cup_{\Sigma} W$ and $V' \cup_{\Sigma'} W'$ are stably equivalent, so we may assume that $\Sigma$ and $\Sigma'$ are isotopic in $M$ after some number of elementary stabilizations.  In addition, let $\Gamma$ be a spine of $V$ such that $\Gamma$ intersects each trivial arc of $\A$ exactly once, so that $\Sigma$ is isotopic to $\pd \ceta(\Gamma)$.  We can find a similar spine $\Gamma'$ for $W'$ intersecting each arc of $\n'$ exactly once, and after a small isotopy we may assume that $\Gamma$ and $\Gamma'$ are disjoint; thus, after isotopy $V$ and $W'$ are disjoint. \\

Now, let $X = \overline{V' - V} = \overline{W - W'}$.  Since $\Sigma$ and $\Sigma'$ are isotopic in $M$, we have that $X = \Sigma \X I = \Sigma' \X I$.  Let $\gamma = K \cap X$ and $\Sigma^* = \Sigma \X \{\frac{1}{2}\}$.  After an isotopy of $\Sigma^*$, we may assume that $\Sigma^*$ is a bridge surface for $(X,\gamma)$.  It is not difficult to see that $\Sigma^*$ is also a bridge surface for $(M,K)$.  We claim that $\Sigma^*$ is the result of elementary perturbations of $\Sigma$ and $\Sigma'$.  Observe that $\Sigma^*$ is a bridge surface for $(V',\A')$.  By Theorem \ref{stand}, $\Sigma^*$ is the result of elementary perturbations of the trivial splitting surface, which is isotopic to $\Sigma'$.  Similarly, $\Sigma^*$ is the result of elementary perturbations of $\Sigma$, completing the proof.

\end{proof}
\end{theorem}

\section{The Dual Curve Complex and the Pants Complex}\label{pnts}

Next, we turn to the pants complex, first defined by Hatcher and Thurston in \cite{hatcher}, and the related dual curve complex, due to Johnson \cite{johnson}.  The curve complex was defined by Harvey in \cite{harvey} and first used to study Heegaard splittings by Hempel in \cite{hempel}.  Let $S$ be a closed surface with genus $g$ and $c$ boundary components.  If we wish to emphasize $g$ and $c$, we will denote $S$ by $S_{g,c}$.  The \emph{curve complex} $C(S)$ of $S$ is a simplicial flag complex whose vertices correspond to isotopy classes of essential simple closed curves in $S$ and whose $k$-cells correspond to collections of $k+1$ isotopy classes of pairwise disjoint essential simple closed curves. \\

\begin{figure}[h]
  \centering
    \includegraphics[width=.9\textwidth]{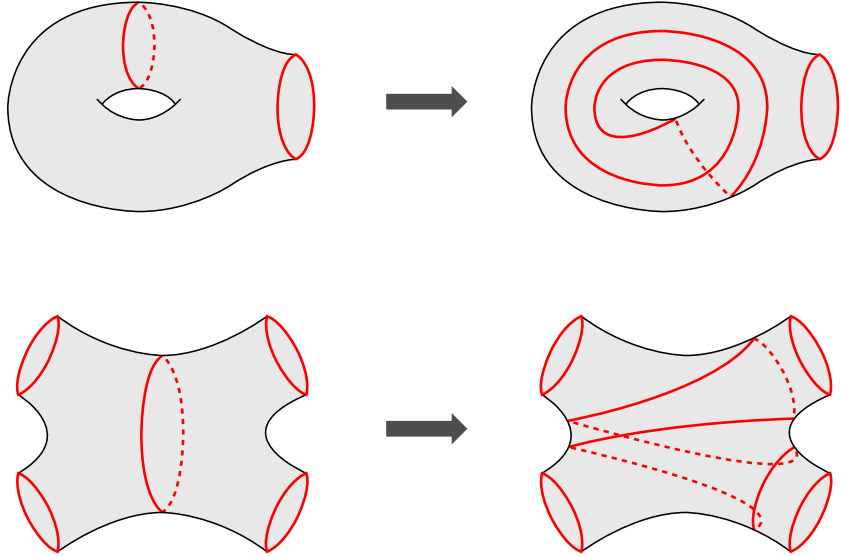}
    \caption{Examples of edges in $C^*(S)$} \label{dual}
\end{figure}

A maximal collection of pairwise disjoint and non-isotopic simple closed curves in a surface $S_{g,c}$ consists of $3g+c-3$ elements; thus, the dimension of $C(S)$ is $3g+c-4$.  The \emph{dual curve complex} $C^*(S)$ is a graph whose vertices correspond to the $3g+c-4$-cells of $C(S)$, where two vertices $v$ and $v'$ are connected by edge if the corresponding $3g+c-4$ cells intersect in a face of codimension one.  In terms of curves in $S$, each vertex $v \in C^*(S)$ represents a collection of $3g+c-3$ isotopy classes of disjoint simple closed curves, which cut $S$ into thrice-punctured spheres (or pairs of pants).  Two vertices share an edge if they differ by a single curve.  Figure \ref{dual} depicts edges in this complex. \\

\begin{figure}[h]
  \centering
    \includegraphics[width=.9\textwidth]{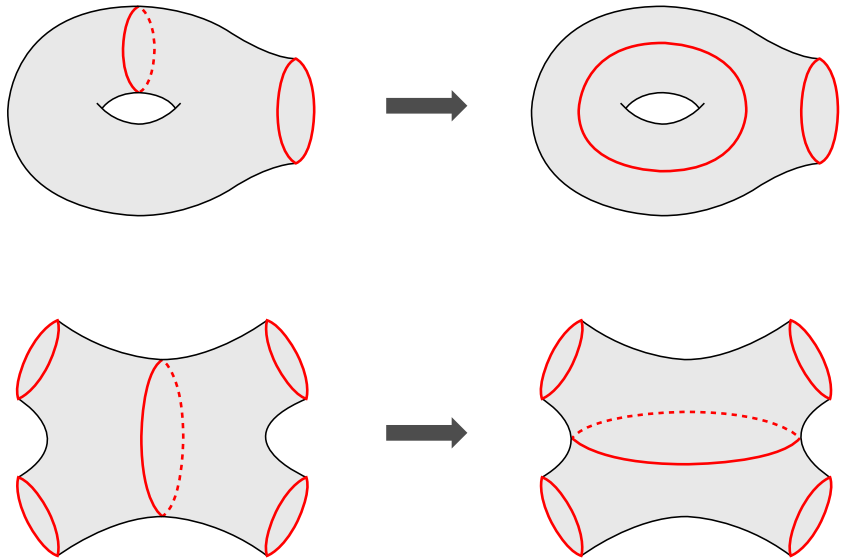}
    \caption{Examples of edges in $P(S)$} \label{pants}
\end{figure}

Each codimension one face $\hat{v}$ of a pants decomposition $v$ can be identified with the curve $u = v - \hat{v}$.  Notice that $S \setminus \hat{v}$ consists of some number of pairs of pants and one component containing $u$ which is either a 4-punctured sphere $S_{0,4}$ or a punctured torus, $S_{1,1}$.  The \emph{pants complex} $P(S)$ is a graph with the same vertex set as $C^*(S)$, where two vertices $v$ and $v'$ are connected by an edge if their corresponding $3g+c-4$-cells share a face of codimension one, where the removed curves $u$ and $u'$ intersect minimally.  That is, if $S \setminus \hat{v}$ contains an $S_{1,1}$ component, we require $u$ and $u'$ to intersect once; if $S \setminus \hat{v}$ has an $S_{0,4}$ component, $u$ and $u'$ must intersect twice.  See Figure \ref{pants}.  Assigning length one to each edge, we construct natural metrics $D(\cdot,\cdot)$ and $D^P(\cdot,\cdot)$ on the vertices of the dual curve complex and pants complex, respectively. \\

In general, we will let $v$ refer to a vertex in $C^*(S)$ and to a collection of $3g+c - 3$ simple closed curves, depending on the context. \\

\section{The distance of a bridge splitting}\label{dist}

For the remainder of the paper, we suppose that $M$ is a closed, orientable 3-manifold.  We do not assume that $M$ is irreducible, but we suppose that $M$ does not contain a non-separating embedded $S^2$ (equivalently, the prime decomposition of $M$ contains no $S^2 \X S^1$ summands).  The concepts presented may be extended to manifolds with boundary or non-seperating 2-spheres, but we omit this extension for simplicity.  For further reference, see \cite{johnson}. \\    

We will use the dual curve and pants complexes to define the distance of a $(g,b)$-bridge splitting given by $(M,K) = (V,\A) \cup_{\Sigma} (W,\n)$.  A \emph{compressing disk} in a 3-manifold $N$ with boundary is a properly embedded disk $D$ such that $\pd D$ is an essential simple closed curve in $\pd N$.  A \emph{compressing disk} in $(V,\A)$ (or $(W,\n)$) is a compressing disk $D$ in $V - \eta(\A)$ such that $\pd D \subset \pd V$.  Lastly, a \emph{cut disk} in $(V,\A)$ (or $(W,\n)$) is a properly embedded disk $D$ such that $\text{int}(D) \cap \A$ is a single point and $\pd D$ is essential in $\pd V- \eta(\pd A)$. \\

Observe that $\Sigma$ may be interpreted as a genus $g$ surface with $2b$ punctures $S_{g,2b}$ (which we will denote $\Sigma_K$); hence, a pants decomposition of $\Sigma_K$ will contain $3g + 2b -3$ essential simple closed curves.  We say that a vertex $v \in C^*(\Sigma_K)$ \emph{defines} $(V,\A)$ (or $(W,\n)$) if every curve of $v$ bounds a compressing disk or a cut disk in $(V,\A)$.  Then, the \emph{distance} $D(\Sigma)$ of the bridge splitting $(M,K) = (V,\A) \cup_{\Sigma} (W,\n)$ is defined to be the minimum distance $D(v,w)$ taken over all pairs of vertices $v,w \in C^*(\Sigma_K)$ defining $(V,\A)$ and $(W,\n)$, respectively.  Similarly, the \emph{pants distance} $D^P(\Sigma)$ is the minimum distance $D^P(v,w)$ between vertices $v,w \in P(\Sigma_K)$ defining $(V,\A)$ and $(W,\n)$, respectively.  These definitions resemble the definitions of the distance and pants distance of a Heegaard splitting of a manifold introduced by Johnson in \cite{johnson}. \\

Since every edge in $P(\Sigma_K)$ is also an edge in $C^*(\Sigma_K)$, we see immediately that
\[ D(\Sigma) \leq D^P(\Sigma).\]
In general, if $(M,K) = (V,\A) \cup_{\Sigma} (W,\n)$ is a bridge splitting, and $v,w$ are elements of $C^*(\Sigma_K)$ such that $v$ defines $(V,\A)$, $w$ defines $(W,\n)$, and $D(v,w) = D(\Sigma)$, we say that the pair $(v,w)$ \emph{yields} $D(\Sigma)$.  If $(v,w)$ yields $D(\Sigma)$, every curve in $v \cap w$ bounds compressing disks or cut disks in both $(V,\A)$ and $(W,\n)$; hence, this collection of curves cannot be too large.  Specifically, if $\Sigma$ gives a $(g,b)$-splitting, then $v$ and $w$ must differ by at least $g+b-1$ curves.  This is made clear by the next lemma.

\begin{lemma}\label{L1}
Suppose that $(M,K) = (V,\A) \cup_{\Sigma} (W,\n)$ is a $(g,b)$-bridge splitting.  Then
\[ D(\Sigma) \geq g + b - 1.\]
\begin{proof}
Suppose that $(v,w)$ yields $D(\Sigma)$, and let $\mathcal{N}$ denote the collection of non-separating curves in $v$.  Since $M$ contains no non-separating 2-sphere, $\mathcal{N} \cap w =\emp$; hence $D(v,w) \geq |\mathcal{N}|$.  We claim that $|\mathcal{N}| \geq g$.  First, let $g \geq 2$.  Since $v$ decomposes $\Sigma_K$ into pairs of pants, $v$ must decompose the closed surface $\Sigma$ into pairs of pants, annuli, and disks.  In particular, some subcollection $v' \subset v$ is a pants decomposition of $\Sigma$.  By an Euler characteristic argument, at least $g$ curves in $v'$ are nonseparating.  If $g = 1$, $v$ must contain at least one nonseparating curve.  In any case, $|\mathcal{N}| \geq g$. \\

Now if $g=0$ and $b = 2$, both $v$ and $w$ consist of exactly one curve.  If $v = w$, then this curve is the intersection of $\Sigma$ with an essential 2-sphere or an essential annulus in $E(K)$, neither of which can occur for a 2-bridge knot $K$; thus, $D(\Sigma) \geq 1$. \\

We prove the lemma by induction on the ordered pair $(b,g)$, with the dictionary ordering.  The cases for which $b=1$ and $(b,g) = (2,0)$ have been completed above.  Thus, suppose that the lemma is true for any $(g_0,b_0)$-bridge splitting of an arbitrary $(M_0,K_0)$, where $(b_0,g_0) < (b,g)$.  Let $\mathcal{S} = v \cap w$. \\

If any curve $\gamma$ in $\mathcal{S}$ bounds cut disks in $(V,\A)$ and $(W,\n)$, then there exists a decomposing sphere $S$ such that $|S \cap K| = 2$ and $S \cap \Sigma = \gamma$.  Cutting $(M,K)$ along $S$ induces two new splittings: $(g_1,b_1)$- and $(g_2,b_2)$-bridge splittings of $(M_1,K_1)$ and $(M_2,K_2)$ with bridge surfaces $\Sigma_1$ and $\Sigma_2$, where $g = g_1 + g_2$ and $b = b_1 + b_2 - 1$.  Suppose without loss of generality that $b_1 \leq b_2$.  If $b_1 = 1$, then $g_1 > 0$ so $g_2 < g$.  If $b_1 > 1$, then $b_2 < b$.  In either case, $(b_1,g_1) < (b,g)$ and $(b_2,g_2) < (b,g)$; thus it follows by induction that
\[ D(\Sigma) = D(\Sigma_1) + D(\Sigma_2) \geq g_1 + b_1 - 1 + g_2 + b_2 - 1 = g + b - 1.\]

Thus, suppose that every curve in $\mathcal{S}$ bounds compressing disks in $V$ and $W$.  If $|\mathcal{S}| \leq 2g - 3$, then $v$ and $w$ differ by at least $g +2b$ curves; hence $D(\Sigma) \geq g + 2b > g + b - 1$.  On the other hand, if $|\mathcal{S}| > 2g - 3$, then $|\mathcal{N} \cup \mathcal{S}| > 3g - 3$, so either there exists there exists $\gamma \in \mathcal{S}$ that bounds a disk $D$ in $\Sigma$ or there exist $\gamma_1,\gamma_2 \in \mathcal{S}$ that are isotopic in $\Sigma$.  In the first case, $D$ contains each of the $2b$ punctures of $\Sigma_K$, which means that disjoint collections of $2b - 2$ curves in $v$ and $2b-2$ curves in $w$ are contained in $\text{int}(D)$, and none of these curves are in $\mathcal{S} \cup \mathcal{N}$.  It follows that
\[ D(\Sigma) \geq g + 2b - 2 > g + b - 1,\]
since $b \geq 2$. \\

In the second case, suppose that no curve in $\mathcal{S}$ bounds a disk in $\Sigma$, so that there exist $\gamma_1,\gamma_2 \in \mathcal{S}$ which cobound an annular component $A$ of $\Sigma \setminus\mathcal{S}$.  This annulus must contain each of the $2b$ punctures; hence disjoint collections of $2b - 1$ curves in $v$ and $2b-1$ curves in $w$ are contained in $\text{int}(A)$, and
\[ D(\Sigma) \geq g + 2b- 1 >  g +b -1,\]
completing the proof of the lemma.
\end{proof}
\end{lemma}

We can also use the above arguments to prove the following:

\begin{lemma}\label{L2}
Suppose that $(M,K) = (V,\A) \cup_{\Sigma} (W,\n)$ is a $(g,b)$-bridge splitting.  If
\[ D(\Sigma) = g + b - 1,\]
then $M$ is $S^3$ or a connected sum of lens spaces and $K$ is the unknot or a connected sum of 2-bridge knots contained in a 3-ball.
\begin{proof}
Suppose $(v,w)$ yields $D(\Sigma)$.  As above, let $\mathcal{S} = v \cap w$.  The statement of the lemma is trivial for $(g,b) = (0,2)$.  Let $(g,b) = (1,1)$, so that each of $v$ and $w$ contain two curves.  In this case, $M$ is clearly $S^3$ or a lens space.  At least one curve in $v$, call it $v_1$, is non-separating; hence $v_1 \notin w$.  It follows that the other curve $v_2$ is in $v \cap w$, and so $v_2$ must bound compressing disks $D_1$ in $(V,\A)$ and $D_2$ in $(W,\n)$.  Now $S = D_1 \cup D_2$ bounds a ball containing $K$ in a 1-bridge position, implying that $K$ is the unknot. \\

Next, suppose that $b =1$ and $g \geq 2$.  By an argument in Lemma \ref{L1}, a subcollection $\mathcal{N}_v \subset v$ of $g$ curves bound non-separating compressing disks in $V$, and a disjoint set $\mathcal{N}_w \subset w$ of $g$ curves bound non-separating disks in $W$.  By assumption, $D(\Sigma) = g$, so $v = \mathcal{N}_v \cup \mathcal{S}$ and $w = \mathcal{N}_w \cup \mathcal{S}$.  It follows from Lemma 11 of \cite{johnson} that $M$ is $S^3$ or a connected sum of lens spaces. \\

If no curve of $\mathcal{S}$ bounds a cut disk in $(V,\A)$, then some curve $\gamma \in \mathcal{S}$ cobounds a pair of pants with the two punctures of $\Sigma_K$.  This implies that $\gamma$ is the intersection of an essential 2-sphere $S$ in $E(K)$ with $\Sigma$, and $S$ bounds a ball in $M$ which contains $K$.  Thus $K$ is a 1-bridge knot in a 3-ball, so $K$ is the unknot.  On the other hand, if some curve $\gamma \in \mathcal{S}$ bounds cut disks in $(V,\A)$ and $(W,\n)$, then cutting along $\gamma$ induces $(g_1,1)$- and $(g_2,1)$-bridge splittings of $(M_1,K_1)$ and $(M_2,K_2)$ with bridge surfaces $\Sigma_1$ and $\Sigma_2$, where $g = g_1 + g_2$ and $M = M_1 \# M_2$.  The statement follows by induction on $g$. \\

Finally, as above we induct on $(b,g)$ with the dictionary ordering.  The base cases $(g,b) = (0,2),(g,1)$ have been shown.  Suppose that $b \geq 2$, and if $b = 2$, $g \geq 1$.  If no curve in $\mathcal{S}$ bounds cut disks in $(V,\A)$ and $(W,\n)$, then $D(\Sigma) > g + b  -1$ by Lemma \ref{L1}, a contradiction.  It follows that some curve $\gamma \in \mathcal{S}$ bounds cut disks in $(V,\A)$ and $(W,\n)$, so that cutting along $\gamma$ induces $(g_1,b_1)$- and $(g_2,b_2)$-bridge splittings of $(M_1,K_1)$ and $(M_2,K_2)$.  Using the inductive hypothesis as in Lemma \ref{L1}, the lemma follows.
\end{proof}
\end{lemma}

Similar statements can be shown for $D^P(\Sigma)$:

\begin{lemma}\label{L3}
Suppose that $(M,K) = (V,\A) \cup_{\Sigma} (W,\n)$ is a $(g,b)$-bridge splitting.  Then
\[ D^P(\Sigma) \geq g + b - 1.\]
\begin{proof}
This follows from Lemma \ref{L1} and the fact that $D^P(\Sigma) \geq D(\Sigma)$.
\end{proof}
\end{lemma}

\begin{lemma}\label{L4}
Suppose that $(M,K) = (V,\A) \cup_{\Sigma} (W,\n)$ is a $(g,b)$-bridge splitting.  If
\[ D^P(\Sigma) = g + b - 1,\]
then $M$ is $S^3$ and $K$ is the unknot.
\begin{proof}
This proof is essentially the same as that of Lemma \ref{L2}, and so we omit most of the details.  Suppose that $(v,w)$ yields $D^P(\Sigma)$, and let $\mathcal{S} = v \cap w$.  If $(g,b) = (0,2)$, then $v$ is a single curve $v_1$ and $w$ is a single curve $w_1$, where $|v_1 \cap w_1| = 2$.  In this case, there does not exist a cut disk for $(V,\A)$ or $(W,\n)$, so both $v_1$ and $w_1$ bound compressing disks.  Let $\A = \A_1 \cup \A_2$ and $\n = \n_1 \cup \n_2$.  Then $\A_1$ cobounds a bridge disk with an arc $\A^*_1 \subset \Sigma$, and $v_1$ is isotopic to $\pd \ceta(\A^*_1) \cap \Sigma$.  Similarly, $\n_1$ cobounds a bridge disk with an arc $\n^*_1 \subset \Sigma$, and $w_1$ is isotopic to $\pd \ceta(\n^*_1) \cap \Sigma$.  Since $|v_1 \cap w_1| = 2$, we have $|\A^*_1 \cap \n^*_1| = 1$, so $\Sigma$ is perturbed and $K$ is the unknot. \\

In the case $(g,b) = (1,1)$, we have that $v$ contains a meridian $v_1$ of $V$ and $w$ contains a meridian $w_1$ of $W$, where $|v_1 \cap w_1| = 1$, so $M = S^3$.  In the case $g \geq 2$ and $b = 1$, we can use Lemma 12 from \cite{johnson} instead of Lemma 11 as in the proof of Lemma \ref{L2} to show that $M = S^3$.  For the remainder of the proof, refer to Lemma \ref{L2}.

\end{proof}
\end{lemma} 

\section{The Bridge Complexity of $K$}\label{BC}

\begin{figure}[h]
  \centering
    \includegraphics[width=1.0\textwidth]{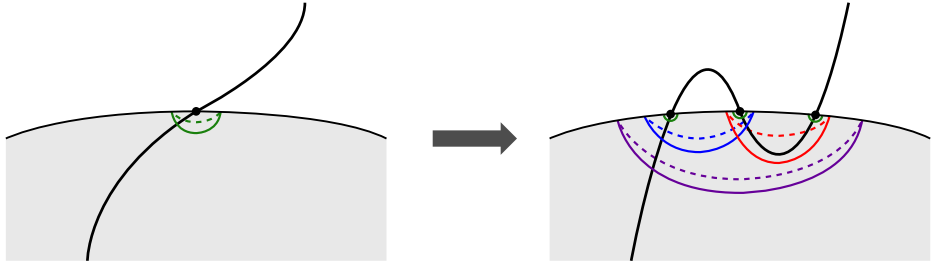}
    \caption{Perturbing and augmenting a vertex in $C^*(\Sigma_K)$ near a puncture, where punctures are green, curves that bound disks on one side are red, curves that bound disks on the other side are blue, and curves that bound disks on both sides are purple} \label{pert2}
\end{figure}

As in \cite{johnson}, we can normalize the distance of any bridge splitting by subtracting a term related to the genus and bridge number of the splitting.  This allows us to define a complexity that is bounded below but does not increase under elementary perturbations or elementary stabilizations.  As above, fix a $(g,b)$-bridge splitting $(M,K) = (V,\A) \cup_{\Sigma} (W,\n)$.  Define the \emph{bridge complexity} and \emph{pants complexity} of $\Sigma$ by
\[ B(\Sigma) = D(\Sigma) - g - b + 1 \qquad \text{ and } \qquad B^P(\Sigma) = D^P(\Sigma) - g - b + 1,\]
respectively.  By Lemmas \ref{L1} and \ref{L3}, we have $B(\Sigma) \geq 0$ and $B^P(\Sigma) \geq 0$.  For $h \geq g$ and $c \geq b$, let $\Sigma^h_c$ be an $(h,c)$-stabilization of $\Sigma$, with the corresponding $(h,c)$-bridge splitting given by $(M,K) = (V^h,\A_c) \cup_{\Sigma^h_c} (W^h,\n_c)$.
\begin{lemma}\label{L5}
For all $h \geq g$ and $c \geq b$, $B(\Sigma^h_{c+1}) \leq B(\Sigma^h_c)$.
\begin{proof}
Suppose that $(v,w)$ yields $D(\Sigma^h_c)$, with $v = v_0,v_1,\dots,v_n = w$ a path of minimal length in $C^*((\Sigma^h_c)_K)$.  Fix a point $x$ of $K \cap \Sigma^h_c$, and let $S = \ceta(x)$ and $s = \pd S$, so that $s$ is a 2-sphere that intersects $K$ in two points and $\Sigma^h_c$ in a single simple closed curve.  Let $\gamma = K \cap S$, and replace $\gamma$ with an arc $\gamma'$ such that $\pd \gamma = \pd \gamma'$, $\gamma \sim \gamma'$, and $\gamma'$ intersects $\Sigma^h_c$ three times.  Then the image of $\Sigma^h_c$ under an isotopy taking $\gamma$ to $\gamma'$ is $\Sigma^h_{c+1}$, the result of performing an elementary perturbation on $\Sigma^h_c$. \\

Note that $S$ contains two bridge disks with respect to the splitting surface $\Sigma^h_{c+1}$, $\Delta \subset V^h$ and $\Delta' \subset W^h$.  Let $\delta = \Sigma^h_{c+1} \cap \pd (\ceta(\Delta))$ and $\delta' = \Sigma^h_{c+1} \cap \pd (\ceta(\Delta'))$, so that $\delta$ bounds a compressing disk in $(V^h,\A_{c+1})$, $\delta'$ bounds a compressing disk in $(W^h,\n_{c+1})$, and $|\delta \cap \delta'| = 2$.  Let $\sigma = s \cap \Sigma^h_{c+1}$, and for each pants decomposition $v_i$, construct a decomposition of $\Sigma^h_{c+1}$ by letting
\[ v_i' = v_i \cup \{\sigma,\delta\}.\]
A depiction of this process appears in Figure \ref{pert2}.  It is clear that $\sigma$ bounds cut disks in both $(V^h,\A_{c+1})$ and $(W^h,\A_{c+1})$, so $v' = v_0'$ defines $(V^h,\A_{c+1})$.  Further, $v' = v_0',v_1',\dots,v_n'$ is a path in $C^*((\Sigma^h_c)_K)$.  Let $v_{n+1}' = v_n' - \{\delta\} \cup \{\delta'\}$.  Then $v_{n+1}'$ defines $(W^h,\n_{c+1})$, and we have $D(\Sigma^h_{c+1}) \leq D(\Sigma^h_c) + 1$.  It follows that
\[ B(\Sigma^h_{c+1}) \leq B(\Sigma^h_c).\]
\end{proof}
\end{lemma}

A similar statement can be made for $B^P$:
\begin{lemma}\label{L6}
For all $c \geq b$, $B^P(\Sigma^h_{c+1}) \leq B^P(\Sigma^h_c)$.
\begin{proof}
The proof is the same as the proof of Lemma \ref{L5}, noting that $|\delta \cap \delta'| = 2$.
\end{proof}
\end{lemma}
As a consequence, we see that neither complexity can increase under perturbation; thus for fixed $h$ both sequences $\{B(\Sigma^h_c)\}$ and $\{B^P(\Sigma^h_c)\}$ converge as $c \rightarrow \infty$.  The situation is similar for stabilizations.

\begin{figure}[h]
  \centering
    \includegraphics[width=1.0\textwidth]{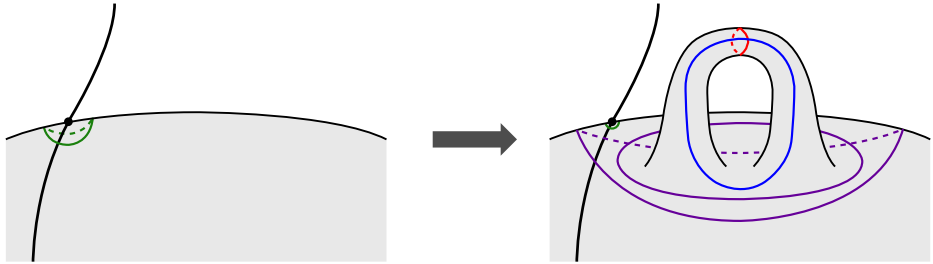}
    \caption{Stabilizing and augmenting a vertex in $C^*(\Sigma_K)$ near a puncture, with coloring as in Figure \ref{pert2}} \label{stab2}
\end{figure}

\begin{lemma}\label{L7}
For all $h \geq g$ and $c \geq b$, $B(\Sigma^{h+1}_c) \leq B(\Sigma^h_c)$.
\begin{proof}
Suppose that $(v,w)$ yields $D(\Sigma^h_c)$, with $v = v_0,v_1,\dots,v_n = w$ a minimal path in $C^*((\Sigma^h_c)_K)$.  Fix a point $x \in K \cap \Sigma^h_c$, let $S = \ceta(x)$, and let $\sigma = \pd(\Sigma^h_c \cap S)$.  Let $\gamma$ be a trivial arc in $S \cap V^h$ bounding a bridge disk $\Delta$ that misses $\A_c$.  Then defining $V^{h+1} = V^h- \eta(\gamma)$ and $W^{h+1} = W \cup \ceta(\gamma)$, we have $(V^{h+1},\A_c) \cup (W^{h+1},\n_c)$ is a stabilization of the splitting given by $\Sigma^h_c$ and has splitting surface $\Sigma^{h+1}_c$.  Note that $\sigma$ bounds cut disks in both $(V^{h+1},\A_c)$ and $(W^{h+1},\n_c)$. \\

We construct a path between $v'$ defining $(V^{h+1},\A_c)$ and $w'$ defining $(W^{h+1},$ $\n_c)$ as follows: define $\tau = \pd (\ceta(\Delta) \cap \Sigma^{h+1}_c)$ and $\eps = \Delta \cap \Sigma^{h+1}_c$, so that $\tau$ bounds compressing disks in both $(V^{h+1},\A_c)$ and $(W^{h+1},\n_c)$ and $\eps$ bounds a compressing disk in $(V^{h+1},\A_c)$.  Let $\eps'$ denote a meridian curve of $\ceta(\gamma)$, so that $\eps'$ bounds a compressing disk in $(W^{h+1},\n_c)$, and $|\eps \cap \eps'| = 1$.  For $0 \leq i \leq n$, let
\[ v_i' = v_i \cup \{\sigma,\tau,\eps\},\]
and let
\[v_{n+1}' = v_n' - \{\eps\} \cup \{\eps'\}.\]
An image of this process can be seen in Figure \ref{stab2}.  Then $v' = v_0'$ defines $(V^{h+1},\A_c)$, $w' = v_{n+1}'$ defines $(W^{h+1},\n_c)$, and $D(v',w') \leq n+1$.  We conclude that
\[B(\Sigma^{h+1}_c) \leq B(\Sigma^h_c),\]
as desired.

\end{proof}
\end{lemma}

\begin{lemma}\label{L8}
For all $h \geq g$ and $c \geq b$, $B^P(\Sigma^{h+1}_c) \leq B^P(\Sigma^h_c)$.
\begin{proof}
The proof is the same as that of Lemma \ref{L6}, noting that $|\eps \cap \eps'|= 1$.
\end{proof}
\end{lemma}
It follows that for fixed $c$ the sequences $\{B(\Sigma^h_c)\}$ and $\{B^P(\Sigma^h_c)\}$ converge as $h \rightarrow \infty$.  We will show convergence as $h,c \rightarrow \infty$ using the next lemma.

\begin{lemma}\label{L9}
There exist $h_* \geq g$ and $c_* \geq b$ such that for every $h \geq h_*$ and $c \geq c_*$,
\[ B(\Sigma^h_c) = B(\Sigma^{h_*}_{c_*}).\]
A similar statement holds for $B^P$.
\begin{proof}
Using Lemmas \ref{L5} and \ref{L7}, we see that the sequence $\{B(\Sigma^{g+n}_{b+n})\}_{n=0}^{\infty}$ is nonincreasing, and by Lemma \ref{L4}, it is bounded below.  Thus, there exists $n_*$ so that $B(\Sigma^{g+n}_{b+n}) = B(\Sigma^{g+n_*}_{b+n_*})$ whenever $n \geq n_*$. \\

Define $h_* = g + n_*$ and $c_* = b + n_*$, and let $h \geq h_*$ and $c \geq c_*$.  If $h - g = c - b$, then there exists $n \geq n_*$ such that $h = g+n$ and $c = b+n$.  Hence, $B(\Sigma^h_c) = B(\Sigma^{h_*}_{c_*})$ by the argument above.  Suppose that $h - g > c - b$, let $n_1 = c - b$, and let $n_2 = h-g$, noting that $n_1,n_2 \geq n_*$.  Then $B(\Sigma^{h_*}_{c_*}) = B(\Sigma^{g+n_1}_c)$ by the argument above, and $\Sigma^h_c$ is a stabilization of $\Sigma^{g+n_1}_c$, so
\[ B(\Sigma^h_c) \leq B(\Sigma^{g+n_1}_c)\]
by Lemma \ref{L7}.  In addition, $B(\Sigma^h_{b+n_2}) = B(\Sigma^{h_*}_{c_*})$ and $\Sigma^h_{b+n_2}$ is a perturbation of $\Sigma^h_c$; thus by Lemma \ref{L5}
\[ B(\Sigma^h_{b+n_2}) \leq B(\Sigma^h_c).\]
Together, these inequalities yield $B(\Sigma^{h_*}_{c_*}) \leq B(\Sigma^h_c) \leq B(\Sigma^{h_*}_{c_*})$, and so the two complexities are equal, as desired.  A similar argument holds when $h-g < c-b$. \\

Employing Lemmas \ref{L6} and \ref{L8} proves the lemma for $B^P$.
\end{proof}
\end{lemma}

At last, we arrive at the main theorems.

\begin{theorem}\label{main}
The limits below exist, and
\[ \lim_{h \rightarrow \infty} \left(\lim_{c \rightarrow \infty} B(\Sigma^h_c) \right) = \lim_{c \rightarrow \infty} \left(\lim_{h \rightarrow \infty} B(\Sigma^h_c) \right).\]
Moreover, these limits do not depend on the bridge surface $\Sigma$ and thus define an invariant of $(M,K)$,
\[ B(M,K) = \lim_{h,c \rightarrow \infty} B(\Sigma^h_k).\]
\begin{proof}
The first statement of the theorem follows directly from Lemma \ref{L9}.  By Theorem \ref{stab}, any other bridge surface $\Sigma'$ has a stabilization equivalent to $\Sigma^h_c$ for some $h$ and $c$; thus the limiting value $B(\Sigma^{h_*}_{c_*})$ does not depend on $\Sigma$ and is an invariant of $(M,K)$.

\end{proof}
\end{theorem}

We call $B(M,K)$ the \emph{bridge complexity} of $(M,K)$, where we omit $M$ when $K \subset S^3$.  Of course, we have proved the same statement for $B^P$:
\begin{theorem}
The quantity
\[ B^P(M,K) = \lim_{h,c \rightarrow \infty} B^P(\Sigma^h_c)\]
defines an invariant of $(M,K)$.
\end{theorem}
We call $B^P(M,K)$ the \emph{pants complexity} of $(M,K)$.

\section{Properties of $B$ and $B^P$}\label{prop}

Here we outline some basic properties of the invariants $B$ and $B^P$.  First, we define the invariants presented in \cite{johnson}, the \emph{Heegaard complexity} and \emph{pants complexity} of a 3-manifold $M$.  These notions will seem familiar given the material presented in Section \ref{BC}.  Suppose that $M = V \cup_{\Sigma} W$ is a Heegaard splitting with Heegaard surface $\Sigma$, where the genus of $\Sigma$ is $g$.  A vertex $v \in C^*(\Sigma)$ defines $V$ if every curve bounds a compressing disk, and the distance $D(\Sigma)$ is the length of a shortest path in $C^*(\Sigma)$ between vertices $v$ and $w$ defining $V$ and $W$.  For a genus $h$ stabilization of this splitting with surface $\Sigma^h$, let $A(\Sigma^h) = D(\Sigma^h) - h$.  Then $A(M) = \lim_{h \rightarrow \infty} A(\Sigma^h)$ is an invariant of $M$, called the \emph{Heegaard complexity}.  The \emph{pants complexity} $A^P(M)$ is defined similarly.  Note that $D(\Sigma)$ and $D^P(\Sigma)$ here are defined only for surfaces $\Sigma$ of genus at least 2. \\

The behavior of the sequences $\{A(\Sigma^h)\}$ and $\{A^P(\Sigma^h)\}$ is slightly different than that of $\{B(\Sigma^h_c)\}_{h\rightarrow \infty}$ and $\{B^P(\Sigma^h_c)\}_{h \rightarrow \infty}$ since there are no punctures about which to augment a path in $C^*(\Sigma)$ or $P(\Sigma)$; initially, these sequences may increase by some fixed amount, but they quickly become nonincreasing as soon as a minimal path in $C^*(\Sigma)$ or $P(\Sigma)$ fixes a curve.  This behavior is described in the next lemma.
\begin{lemma}\label{L11a}
Suppose, for some Heegaard splitting $M = V \cup_{\Sigma} W$ of genus $g$ and for a path $v = v_0,v_1,\dots,v_n = w$ in $C^*(\Sigma)$ between vertices $v$ and $w$ defining $V$ and $W$, that there exists a simple closed curve $\mu$ in every $v_i$.  Then for any stabilization $\Sigma^h$ of $\Sigma$,
\[ A(\Sigma^h) \leq n - g \qquad \text{ and } \qquad A^P(\Sigma^h) \leq n - g.\]
\begin{proof}
The proof is similar to the proof of Lemma \ref{L7}, and so we give only a sketch.  Suppose $(v,w)$ yields $D(\Sigma)$, with $v = v_1,\dots,v_n = w$ a minimal path in $C^*(\Sigma)$.  Let $C$ denote the compressing disk $\mu$ bounds in $V$, and let $C \X I$ be a collar of $C$ in $V$ so that $C = C \X \{0\}$.  Define $C' = C \X \{1\}$ and $\mu' = \pd C'$.  Let $\gamma$ be a trivial arc in $C \X I$ with $\Delta$ a bridge disk for $\gamma$.  Note that the Heegaard splitting $(V- \eta(\mu)) \cup (W \cup \ceta(\mu)) = V' \cup_{\Sigma'} W'$ is a stabilization of $\Sigma$.  Let $\tau = \pd (\Sigma' \cap \ceta(\Delta))$ and $\eps = \Delta \cap \Sigma'$, with $\eps'$ a meridian of $\ceta(\mu)$.  For $0 \leq i \leq n$, define
\[ v_i' = v_i \cup \{\mu',\tau,\eps\},\]
and let $v_{n+1}' = v_n - \{\eps\} \cup \{\eps'\}$.  Then $v' = v_0'$ defines $V'$, $w' = v_{n+1}'$ defines $W$, and $v' = v_0',\dots,v_{n+1}' = w'$ is a path in $C^*(\Sigma')$.  Thus $D(\Sigma') \leq n + 1$, and the statement of the lemma follows after finitely many repetitions of this process (and noting for $D^P$ that $|\eps \cap \eps'| = 1$).
\end{proof}
\end{lemma}

We may now relate $B(M,K)$ to $A(M)$ and $B^P(M,K)$ to $A^P(M)$:
\begin{lemma}\label{L11}
For a manifold $M$ containing a knot $K$,
\[ B(M,K) \geq A(M) \qquad \text{ and } \qquad B^P(M,K) \geq A^P(M).\]
\begin{proof}
We prove the statement for $B$ and $A$, since the proof for $B^P$ and $A^P$ is identical.  Suppose that $(M,K) = (V,\A) \cup_{\Sigma} (W,\n)$ is a $(g,b)$-bridge splitting and that $(v,w)$ yields $D(\Sigma)$, with $v = v_0,v_1,\dots,v_n = w$ be a minimal path in $C^*(\Sigma_K)$.  Note that a pants decomposition of $\Sigma_K$ contains $3g + 2b -3$ curves.  By Theorem \ref{main}, for large $g$ and $b$, $D(\Sigma) = B(M,K) + g + b -1$, and so by choosing $g$ and $b$ large enough, we may ensure that $D(\Sigma) < 3g + 2b - 3$, guaranteeing the existence of a curve $\gamma$ common to every $v_i$.  If $\gamma$ is inessential in $\Sigma$, it can be made essential after another elementary stabilization, so we assume without loss of generality that $\gamma$ is essential in $\Sigma$. \\

Let $v' \subset v$ be a subcollection containing $\gamma$ and $3g-4$ other curves that give a pants decomposition of the closed surface $\Sigma$.  Since each curve in $v$ bounds a compressing disk or cut disk in $(V,\A)$, each curve in $v'$ must bound a compressing disk in $V$, so $v'$ defines $V$. \\

Now, for each $i$, choose a subcollection $v_i' \subset v_i$ containing $\gamma$ and $3g-4$ other curves decomposing $\Sigma$ into pairs of pants.  By the same reasoning as above, $w' = v_n'$ defines $W$.  We claim that $v' = v_0',v_1',\dots,v_n' = w'$ is a path in $C^*(\Sigma)$.  Since $v_i \cup v_{i+1}$ contains exactly one pair $(u_i,u_{i+1})$ of curves which are not disjoint, $v_i' \cup v_{i+1}'$ can contain at most one such pair.  Thus, there exists a path fixing $\gamma$ between $v'$ and $w'$ in $C^*(\Sigma)$ of length at most $n$.  It follows from Lemma \ref{L11a} that
\[ A(M) = \lim_{h \rightarrow \infty} A(\Sigma^h) \leq n - g,\]
which proves the lemma in the case that $b=1$.  If $g = 0$, $M = S^3$ and $A(M) = A^P(M) = 0$. \\

As in the proof of Lemma \ref{L1}, we induct on $(b,g)$ with the dictionary order.  Suppose that $u_i \in v_i$, $u_i \notin v_{i+1}$, and $u_i$ cobounds a pants component of $\Sigma_K \setminus v_i$ with a puncture of $\Sigma_K$, call it $x$ and some other curve $y$.  Then in the closed surface $\Sigma$, $u_i$ is isotopic to $y$.  Let $u_i'$ denote the unique curve such that $u_i' \in v_{i+1}$ but $u_i' \notin v_i$.  Then $u_i'$ also cobounds a pants component of $\Sigma_K \setminus v_{i+1}$ with $x$ and some other curve $z$, implying that $u_i'$ is isotopic to $z$ in $\Sigma$.  We conclude that the induced pants decompositions $v'_i$ and $v'_{i+1}$ must be identical. \\

Thus, suppose that $v \cap w = \emp$, and let $\mathcal{T} \subset v$ denote the subcollection of curves $u$ which cobound a pair of pants with one of the punctures of $\Sigma_K$.  Since there are $2b$ such punctures, it follows that $|\mathcal{T}| \geq b$.  As $v \cap w = \emp$, each curve $u \in \mathcal{T}$ satisfies $u \in v_i$ but $u \notin v_{i+1}$ for some $i$; hence we see that $v_i' = v'_{i+1}$ for at least $b$ distinct values of $i$.  It follows that the length of the path $v_0', \dots, v_n'$ is at most $n -b$, and so
\[ A(M) \leq n - g - b \leq n - g - b + 1 = B(M,K).\]

On the other hand, suppose that $v \cap w \neq \emp$, so that there is a curve $\delta$ which bounds compressing or cut disks in both $(V,\A)$ and $(W,\n)$.  If $\delta$ bounds cut disks, then cutting along $\Delta$ induces a $(g_1,b_1)$-splitting of some $(M_1,K_1)$ and $(g_2,b_2)$-splitting of some $(M_2,K_2)$, where $(M,K) = (M_1,K_1) \# (M_2,K_2)$ with $g_1 + g_2 = g$ and $b_1 + b_2 - 1 = b$.  Let $\Sigma_1$ and $\Sigma_2$ denote the corresponding splitting surfaces, so that $D(\Sigma) = D(\Sigma_1) + D(\Sigma_2)$. \\

By induction, $A(M_1) \leq B(M_1,K_1)$ and $A(M_2) \leq B(M_2,K_2)$, and by Lemma 19 of \cite{johnson}, $A(M) \leq A(M_1) + A(M_2)$.  It follows that
\begin{eqnarray*}
A(M) &\leq& B(M_1,K_1) + B(M_2,K_2) \\
&\leq& D(\Sigma_1) - g_1 - b_1 + 1 + D(\Sigma_2) - g_2 - b_2 + 1 \\
&\leq& D(\Sigma) - g - b + 1 = B(M,K).
\end{eqnarray*}
On the other hand, if $\delta$ bounds compressing disks, we again use Lemma 19 of \cite{johnson} to produce a similar string of inequalities.
\end{proof}
\end{lemma}
We remark the the first inequality in Lemma \ref{L11} is sharp if and only if $K$ is the unknot or the connected sum of 2-bridge knots contained in a ball, and the second inequality if sharp if and only if $K$ is the unknot, so that the differences $B(M,K) - A(M)$ and $B^P(M,K) - A^P(M)$ may also be useful measures of the complexity of $K$.  The proof of these facts is left as an exercise for the interested reader. \\

In light of Lemma \ref{L11}, we may view $A(M)$ as $B(M,\emp)$ and $A^P(M)$ as $B^P(M,\emp)$, so that $B$ and $B^P$ are extensions of $A$ and $A^P$.  In addition, Lemma \ref{L11} reveals that there are pairs $(M,K)$ of arbitrarily large complexity, since there exist manifolds with large complexity \cite{johnson}. \\

For manifolds $M_1,M_2$ and knots $K_i \subset M_i$, the connected sum $(M_1,K_1) \#$ $(M_2,K_2)$ is constructed by removing balls $B_i \subset M_i$ intersecting $K_i$ in a single trivial arc $\A_i$, and attaching $\pd B_1$ to $\pd B_2$ so that $\pd \A_1$ coincides with $\pd \A_2$.  As with many knot invariants, it is of interest to examine the behavior of $B$ and $B^P$ under the connected sum operation.

\begin{lemma}\label{L12}
For manifolds $M_1,M_2$ and knots $K_i \subset M_i$,
\begin{eqnarray*}
B((M_1,K_1) \# (M_2,K_2)) &\leq& B(M_1,K_1) + B(M_2,K_2) \\
&\text{and} & \\
B^P((M_1,K_1) \# (M_2,K_2)) &\leq& B^P(M_1,K_1) + B_P(M_2,K_2).
\end{eqnarray*}
\begin{proof}
Again, we prove the lemma for $B$, noting that the same proof holds for $B^P$.  For $i=1,2$, suppose that $(M_i,K_i) = (V_i,\A_i) \cup_{\Sigma_i} (W_i,\n_i)$ is a $(g_i,b_i)$-bridge splitting such that $B(\Sigma_i) = B(M_i,K_i)$, and let $(v^i,w^i)$ yield $D(\Sigma_i)$ with $v^i = v^i_0,v^i_1,\dots,v^i_{n_i} = w^i$ a minimal path in $C^*((\Sigma_i)_K)$.  Fix $p_i \in K_i \cap \Sigma_i$, and let $S_i = \ceta(p_i)$.  Let $\sigma_i = \pd (\Sigma_i \cap S_i)$.  Observe that $S_i$ is a ball containing a trivial arc, and we can perform the connected sum operation along $\pd S_i = C^i_V \cup C^i_W$, where $C^i_V = S_i \cap V_i$, $C^i_W = S_i \cap W_i$, and $\pd C^i_V = \pd C^j_W = \sigma_i$. \\

Let $(V,\A) \cup_{\Sigma} (W,\n)$ denote the induced splitting of $(M_1,K_1) \# (M_2,K_2)$, where $\A = \A_1 \cup \A_2$, $\n = \n_1\cup \n_2$, $V$ is $V_1$ glued to $V_2$ along $C^1_V$ and $C^2_V$, and $W$ is $W_1$ glued to $W_2$ along $C^1_W$ and $C^2_W$.  Note that $\Sigma$ is a $(g_1 + g_2,b_1+b_2-1)$-splitting surface.  For $0 \leq j \leq n_1$, define
\[ v_j' = v^1_j \cup v^2 \cup \{\sigma\},\]
and for $n_1+1 \leq j \leq n_1 + n_2$, define
\[ v_j' = w^1 \cup v^2_j \cup \{\sigma\}.\]
Since $\sigma$ bounds cut disks $C^1_V = C^2_V$ and $C^1_W = C^2_W$ in $V$ and $W$, respectively, it is clear that $v' = v_0'$ defines $(V,\A)$ and $v'_{n_1+n_2}$ defines $(W,\n)$.  It follows that
\begin{eqnarray*}
B((M_1,K_1) \# (M_2,K_2)) &\leq& B(\Sigma) \\
&\leq& n_1+n_2 - (g_1+g_2) - (b_1+b_2-1) +1\\
&=& B(M_1,K_1) + B(M_2,K_2).
\end{eqnarray*}

\end{proof}
\end{lemma}
Similar inequalities holds for $A$ and $A^P$; this is Lemma 19 of \cite{johnson}. \\

We may now undertake our first calculations.
\begin{lemma}\label{L13}
For a manifold $M$ and knot $K \subset M$, $B(M,K) = 0$ if and only if $M = S^3$ or $M$ is a connected sum of lens spaces and $K$ is the unknot or a connected sum of 2-bridge knots contained in a ball.
\begin{proof}
First, suppose that $B(M,K) = 0$.  Then there exists a $(g,b)$-splitting surface $\Sigma$ for $(M,K)$ such that $B(\Sigma) = 0$; that is $D(\Sigma) = g + b - 1$.  By Lemma \ref{L2}, $M = S^3$ or $M$ is a connected sum of lens spaces and $K$ is the unknot or a connected sum of 2-bridge knots contained in a ball. \\

Conversely, suppose first that $M = S^3$ and $K$ is the unknot or a 2-bridge knot, and let $(S^3,K) = (B_1,\A) \cup_{\Sigma} (B_2,\n)$ be a $(0,2)$-bridge splitting of $K$.  Then $\Sigma_K$ is a 4-punctured sphere $S_{0,4}$, so a pants decomposition is a single curve, and every curve is distance one from every other curve in $C^*(\Sigma_K)$.  It follows that $D(\Sigma) = 1$, so $B(\Sigma) = 0$.  Since $B(K) \leq B(\Sigma)$, $B(K) = 0$ as well. \\

Now, if $K$ is the connected sum of 2-bridge knots $K_1,\dots,K_n$, then by a repeated application of Lemma \ref{L12}, $B(K) \leq B(K_1) + \dots + B(K_n) = 0$, so $B(K) = 0$.  Finally, suppose $M$ the lens space $L(p,q)$ and $K$ is the connected sum of 2-bridge knots contained in the ball.  Let $\Sigma$ be a $(g,b)$-splitting surface for $K \subset S^3$, so that $B(\Sigma) = 0$ by the arguments above.  Stabilize $\Sigma$ to $\Sigma'$ as in Lemma \ref{L7} and suppose that $(v,w)$ yields $D(\Sigma')$.  Replacing the newly constructed longitude curve of $w$ with a $(p,q)$-curve gives a path of a the same length and whose endpoints define either side of a splitting of $(M,K)$ with surface $\Sigma''$.  It follows that $B(M,K) \leq B(\Sigma'') = 0$.  A similar construction works when $M$ is an arbitrary sum of lens spaces.
\end{proof}
\end{lemma}
\begin{lemma}\label{L14}
For a manifold $M$ and knot $K \subset M$, $B^P(M,K) = 0$ if and only if $M = S^3$ and $K$ is the unknot.
\begin{proof}

The proof is essentially the same as that of Lemma \ref{L13}, employing Lemma \ref{L4} instead of Lemma \ref{L2} for the forward implication.
\end{proof}
\end{lemma}

\section{Critical components of bridge splittings}\label{crit}

In order to undertake some simple calculations in the next section, we will need to bound the bridge complexity below more sharply than the bound given by Lemma \ref{L11}.  For this reason, we study the pieces of a splitting which are indecomposable (in some sense) and which realize the complexity of $(M,K)$.  Let $M_K$ denote the \emph{exterior} of $K$ in $M$, $M_K = M - \eta(K)$.  If $M = S^3$ we write $E(K) = M_K$.  We say that a pair $(M,K)$ is \emph{irreducible} if two conditions are satisfied:
\be
\item Any embedded 2-sphere in $M_K$ bounds a ball in $M_K$, and
\item Any properly embedded annulus in $M_K$ whose boundary consists of meridians of $K$ is boundary parallel.
\ee
For example, if $K$ is any prime knot in $S^3$, then $(S^3,K)$ is irreducible. \\

Now, suppose that $(M,K)$ is irreducible, and let $(M,K) = (V,\A) \cup_{\Sigma} (W,\n)$ be a $(g,b)$-bridge splitting such that $B(\Sigma) = B(M,K)$ with $(v,w)$ yielding $D(\Sigma)$ and $v = v_0,\dots,v_n = w$ a minimal path in $C^*(\Sigma_K)$.  Further, let $\Sau$ denote the collection of curves $\gamma$ contained in every $v_i$.  Since each $\gamma \in \Sau$ bounds compressing disks or cut disks in $(V,\A)$ and $(W,\n)$, cutting along these disks and cut disks decomposes $(M,K)$ into pieces, some of which -- call them $(M_1,K_1),\dots,(M_l,K_l)$ -- are knots contained in manifolds (these are created by cutting along cut disks), and some of which -- call them $M_{l+1},\dots,M_m$ -- are manifolds (these are created by cutting along compressing disks).  In addition, this process splits $\Sigma$ into bridge surfaces or Heegaard surfaces $\Sigma_1,\dots,\Sigma_m$ for $(M_j,K_j)$ or $M_j$ and possibly pairs of pants bounded by elements of $\Sau$.\\

Let $g_j$ denote the genus and $b_j$ the bridge number of $\Sigma_j$ (in the case that $\Sigma_j$ is a Heegaard surface, set $b_j = 0$).  Then $\sum g_j = g$ and $\sum b_j - l + 1 = b$.  Now, let $e_i$ denote the edge between $v_{i-1}$ and $v_i$ in $C^*(\Sigma_K)$.  For each $i$, there exists exactly one $\Sigma_j$ such that $v_{i-1} \cap \Sigma_j \neq v_i \cap \Sigma_j$.  In this case, we say that the edge $e_i$ is \emph{contained} in $\Sigma_j$.  Now, by the irreducibility of $(M,K)$, we may assume that $(M,K) = (M_1,K_1)$ and $(M_j,K_j) = (S^3,0_1)$, $M_j = S^3$ for $j \geq 2$.  For each $\Sigma_j$, let $n_j$ denote the number of edges contained $\Sigma_j$, so that $\sum n_j = n = D(v,w) = D(\Sigma)$.  Since $\Sigma$ has minimal complexity, each $\Sigma_j$ must also have minimal complexity.  Of course, the complexities of $(S^3,0_1)$ and $S^3$ are zero, which means that for $2 \leq j \leq l$, $n_j = g_j + b_j - 1$ and for $l+1 \leq j \leq m$, $n_j = g_j$.  It follows that
\begin{equation}\label{cc}
B(M,K) = n - g - b + 1 = \sum n_j - \sum g_j - \sum b_j + l  = n_1 - g_1 - b_1 + 1.
\end{equation}
In other words, the splitting surface $\Sigma_1$ for $(M_1,K_1)$ contributes all of the complexity to $B(M,K)$.  The quantity is not a priori the same as $B(\Sigma_1)$, since the process of cutting may leave scars of compressing disks on $\Sigma_1$, in which case $\Sigma_1$ inherits pants decompositions with more punctures and hence more curves than a standard $(g_1,b_1)$-splitting. \\

\begin{figure}[h]
  \centering
    \includegraphics[width=1.0\textwidth]{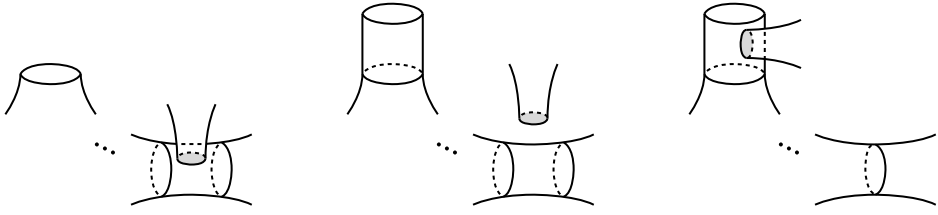}
    \caption{We create $v_i^+$ from $v_i$ by adding a curve $\gamma_0$ parallel to a puncture, gluing a disk to $\gamma$ in $\Sigma_1'$, replacing $\gamma_i$ and $\gamma_i'$ with a single curve, and reattaching $\Sigma_i$ or $P$ along $\gamma$} \label{scar}
\end{figure}

\begin{lemma}\label{L17}
The splitting surface $\Sigma_1$ defined above has no scars from cutting along compressing disks.
\begin{proof}
The main idea of the proof is that if $\Sigma_1$ has scars, a path in the dual curve complex of $\Sigma_1$ may be collapsed to a shorter path without scars.  Recall that in the construction of $\Sigma_1$, a scar corresponds to the remnants of cutting along a compressing disk whose boundary appears in every vertex of a minimal path in the dual curve complex.  Thus, suppose that $\Sigma_1$ has a scar bounded by $\gamma$, let $\Sigma_1'$ denote the surface $\Sigma_1$ punctured by $K$ and scars.  Suppose that $v = v_0,\dots,v_n = w$ is a minimal path between $v$ and $w$, where $n_1$ of the edges in the path are contained in $\Sigma_1$.  After reordering, we may suppose that the edges $e_1,\dots,e_{n_1}$ are those contained in $\Sigma_1$.  As noted above, $B(M,K) = n_1 - g_1 - b_1 + 1$. \\

Observe that $\gamma \in v_i$ for every $i$, and $v_i \cap \Sigma_1$ is a pants decomposition of $\Sigma_1'$.  For each $i$, let $\gamma_i$ and $\gamma_i'$ denote the two curves of $v_i$ that cobound a pair of pants with $\gamma$.  Note that if $\gamma_i$ and $\gamma_i'$ were the same curve for some $i$, it would imply that $\Sigma_1'$ is a once punctured torus $S_{1,1}$, which is impossible as $b_1 \geq 1$.  Thus, $\gamma_i$ and $\gamma_i'$ are always distinct curves. \\

Now, let $\Sigma_1^*$ denote the surface $\Sigma_1'$ with a disk attached to $\gamma$, so that $\Sigma_1^*$ has one fewer puncture that $\Sigma_1'$.  Each pants decomposition $v_i \cap \Sigma_1'$ descends to a pants decomposition $v_i^*$ of $\Sigma_1^*$ by replacing $\gamma_i$ and $\gamma_i'$ (which are parallel in $\Sigma_1'$) with a single curve $\gamma_i^*$.  If $\gamma_{i-1} = \gamma_i$ and $\gamma_{i-1}' = \gamma_i'$, the edge $e_i$ descends to an edge $e_i^*$ between $v_{i-1}^*$ and $v_i^*$.  Otherwise, if $e_i$ corresponds to removing $\gamma_{i-1}$ or $\gamma_{i-1}'$ and replacing it with a different curve, then $v_{i-1}^* = v_i^*$.  Of course, since neither $\gamma_0$ nor $\gamma_0'$ is in $v_n = w$, there is at least one edge of the second type, so if the induced path has $n_1^*$ edges contained in $\Sigma_1$, $n_1^* < n_1$.  \\

We extend the path of $v_i^*$'s to all of $\Sigma$.  Fix a point of intersection $x \in K \cap \Sigma_1$, and let $\gamma_0 = \pd \ceta(x) \cap \Sigma_1$.  Then $\gamma_0$ bounds a disk $D$ in $\Sigma_1$ containing $x$; let $\gamma_0^* \subset \text{int}(D)$ be a curve bounding a disk $D^* \subset D$ such that $x \notin D^*$.  The process of cutting along $\gamma$ leaves two scars, one in $\Sigma_1$, and the other in a different $\Sigma_i$ or pair of pants $P$.  Taking the connected sum of $\Sigma_1$ and $\Sigma_i$ or $P$ along $D^*$ and the scar of $\Sigma_i$ or $P$, for each $i$ we induce a new pants decomposition
\[ v_i^+ = v_i - (v_i \cap \Sigma_1) \cup v_i^* \cup \{\gamma_0,\gamma_0^*\},\]
creating a new path $v^+ = v_0^+,v_1^+,\dots,v_{n^+}^+ = w^+ \in C^*(\Sigma_K)$.  But $v^+$ defines $(V,\A)$ and $w^+$ defines $(W,\n)$, and since $v_i^+ = v_{i+1}^+$ for at least one value of $i$, this path has length less than $n$.  This implies $B(\Sigma) < B(M,K)$, a contradiction.  See Figure \ref{scar}.

\end{proof}
\end{lemma}

The previous lemma implies that $B(M,K) = B(\Sigma) = B(\Sigma_1)$.  We call $\Sigma_1$ the \emph{critical component} of $\Sigma$ with respect to the path from $v$ to $w$, noting that $\Sigma_1$ depends on this choice of path.  We have the following crucial property:
\begin{lemma}\label{L15}
Let $(M,K)$ be irreducible and suppose $\Sigma_1$ is a $(g_1,b_1)$-splitting surface which is the critical component of some bridge splitting of $(M,K)$.  Then
\[ B(M,K) \geq 2g_1 + b_1 - 2.\]
\begin{proof}
Suppose that $\Sigma_1$ is the critical component of $(V,\A) \cup_{\Sigma} (W,\n)$ with respect to the minimal path $v = v_0,\dots,v_n = w$ between the pair of vertices $(v,w)$ yielding $D(\Sigma)$, and $\Sigma_1$ contains $n_1$ edges.  Any pants decomposition of $\Sigma_1$ contains $3g_1 + 2b_1 - 3$ curves. By the construction of $\Sigma_1$, no curve in $v \cap \Sigma_1$ is in $w \cap \Sigma_1$; hence $\Sigma_1$ contains at least $3g_1 + 2b_1 - 3$ edges.  Hence, using the formula (\ref{cc}) for $B(M,K)$ appearing above,
\[ B(M,K) = n_1 - g_1  -b_1 + 1 \geq 3g_1 + 2b_1 - 3 - g_1 - b_1 + 1 = 2g_1 - b_1 -2.\]

\end{proof}
\end{lemma}

Lemma \ref{L15} can be used to give a lower bound for $B(M,K)$ based on the tunnel number of $K$.  The \emph{tunnel number} a knot $K$ in a manifold $M$, denoted $t(K)$, is the minimal number of arcs that may be attached to $K$ so that the complement in $M$ of a neighborhood of the union of $K$ and the arcs is a handlebody.  Equivalently, $t(K) = g(M_K) - 1$, where $g(M_K)$ is the Heegaard genus of $M_K$.  Lemma \ref{L11} implies that there exist knots of arbitrarily high complexity provided that the complexity of the corresponding manifold containing them is high, but as expected, low complexity manifolds such as $S^3$ also contain knots of arbitrarily high complexity. \\

Any $(g,b)$-bridge splitting of $(M,K)$ gives rise to a Heegaard splitting of $M_K$ via the following construction, which is well-known (see for instance, \cite{hempel}):  Suppose $(M,K) = (V,\A) \cup_{\Sigma} (W,\n)$ is a $(g,b)$-bridge splitting.  Let $V' =V - \eta(\A)$.  Since each arc of $\A$ is trivial, $V'$ is a genus $g+b$ handlebody.  Now, let $W^* = W \cup \ceta(\A)$, so that $W^*$ can be viewed as $W$ with $b$ 1-handles attached, and $W^*$ is also a genus $g+b$ handlebody.  Finally, note that $K \subset W^*$ and any meridian of the attached 1-handles intersects $K$ once, so $W' = W^* - \eta(K)$ is a compression body.  Letting $\Sigma' = \pd_+ W' = \pd_+ V'$, we have $M_K = V' \cup_{\Sigma'} W'$ is a genus $g + b$ Heegaard splitting, from which it follows that $t(K) \leq g + b - 1$.  This inequality is important in establishing another lower bound for $B(M,K)$ (compare this to Lemma 23 of \cite{johnson}):

\begin{lemma}
Suppose $(M,K)$ is irreducible.  Then $B(M,K) \geq t(K) - 1$.
\begin{proof}
Let $(M,K) = (V,\A) \cup_{\Sigma} (W,\n)$ be a bridge splitting such that $B(\Sigma) = B(M,K)$, and let $\Sigma_1$ be the critical component of $\Sigma$ with respect to some minimal path between vertices $v,w \in C^*(\Sigma_K)$ defining $(V,\A),(W,\n)$.  Suppose that $\Sigma_1$ is a $(g_1,b_1)$-bridge splitting for $(M,K)$.  By the argument preceding the lemma, $t(K) \leq g_1 + b_1 - 1$.  It follows by Lemma \ref{L15} that
\[ B(M,K) \geq 2g_1 + b_1 - 2 \geq g_1 + b_1 - 2 \geq t(K) - 1.\]
\end{proof}
\end{lemma}
It follows trivially that for such $(M,K)$, $B^P(M,K) \geq t(K) - 1$ as well.  It also follows that manifolds $M$ with small $A(M)$ (respectively $A^P(M)$) contain knots $K$ so that $B(M,K)$ (respectively $B^P(M,K)$) is arbitrarily large. \\

\section{Some computations}\label{calc}

To begin, we employ the lemmas of Section \ref{crit} to find $B(K_{p,q})$ for $K_{p,q}$ a $(p,q)$-torus knot in $S^3$.  We will always suppose without loss of generality that $p < q$.  The \emph{bridge number} $b(K)$ of a knot $K$, defined by Schubert in \cite{schub}, is the minimum $b$ such that $(S^3,K)$ admits a $(0,b)$-bridge splitting.  In \cite{schult}, Schultens shows that $b(K_{p,q}) = p$.  Thus, by Lemma \ref{L13}, $B(K_{2,q}) = 0$, since $K_{2,q}$ is a 2-bridge knot. \\

Suppose that $(M,K) = (V,\A) \cup_{\Sigma} (W,\n)$ is a bridge splitting such that there exist disjoint curves $\gamma_1$ and $\gamma_2$ in $\Sigma$ such that $\gamma_1$ bounds a compressing disk or a cut disk in $(V,\A)$ and $\gamma_2$ bounds a compressing disk or a cut disk in $(W,\n)$.  In this case, the bridge surface $\Sigma$ is called \emph{$c$-weakly reducible}, a definition due to Tomova \cite{tomova2}.  A knot $K \subset M$ is called \emph{meridionally small} if its exterior $M_K$ contains no essential meridional surfaces.  Every meridionally small knot is prime, and so if $K$ is meridionally small, $(S^3,K)$ is irreducible.  A special case of the main theorem of \cite{tomova2} is the following lemma:

\begin{lemma}\label{L18}
Suppose $K$ is meridionally small, $(S^3,K) = (B_1,\A) \cup_{\Sigma} (B_2,\n)$ is a $(0,b)$-bridge splitting, and the bridge surface $\Sigma$ is $c$-weakly reducible.  Then $\Sigma$ is perturbed.
\end{lemma}

Let $K$ be a knot in $S^3$, with $\Sigma$ be a $(0,b)$-bridge surface for some $b$, and define
\[B_0(\Sigma) = \lim_{c \rightarrow \infty} B(\Sigma^0_c)\]
and
\[B_0^P(\Sigma) = \lim_{c \rightarrow \infty} B^P(\Sigma^0_c).\]
By Theorem \ref{stab}, we know that any two classical bridge splittings (with underlying Heegaard surface a ball) have a common perturbation; hence $B_0(\Sigma)$ and $B_0^P(\Sigma)$ do not depend on $\Sigma$ and define invariants of $K$, so we write $B_0(K) = B_0(\Sigma)$ and $B_0^P(K) = B_0^P(\Sigma)$.  It follows immediately from the Lemmas in Section \ref{BC} that $B(K) \leq B_0(K)$ and $B^P(K) \leq B_0^P(K)$.  Lemma \ref{L18} factors into the following calculation:
\begin{lemma}
If $K$ is meridionally small and $b(K) \geq 3$, then $B_0(K) \geq b(K) - 1$.
\begin{proof}\label{L19}
Let $K = (B_1,\A) \cup_{\Sigma} (B_2,\n)$ be a $(0,b)$-bridge splitting such that $B(\Sigma) = B_0(K)$, and let $\Sigma_1$ be the critical component of $\Sigma$ with respect to a minimal path $v = v_0,\dots,v_n = w$ between the pair of vertices $(v,w)$ yielding $D(\Sigma)$, where $\Sigma_1$ is a $(0,b_1)$-bridge surface.  If $b_1 > b(K)$, then $b_1 \geq b(K) + 1$, and so using the proof of Lemma \ref{L15},
\[ B(\Sigma) \geq b_1 - 2 \geq b(K) - 1.\]
Otherwise, suppose that $b_1 = b(K)$.  A pants decomposition of $\Sigma_1$ contains $2b_1 - 3$ curves.  Suppose that $\Sigma_1$ contains $n_1 \leq 4b_1 - 8$ edges, and reorder the path from $v$ to $w$ so that all the edges between $v_0,\dots,v_{n_1}$ are contained in $\Sigma$.  Now $v \cap v_{2b_1-4}$ contains at least one curve $\gamma_1 \subset \Sigma_1$, which necessarily bounds a compressing disk or a cut disk in $(B_1,\A)$.  Noting that $w \cap \Sigma_1 = v_{n_1} \cap \Sigma_1$ and $n_1 \leq 2b_1 - 8$, we have that $w \cap v_{2b_1 - 4}$ contains at least one curve $\gamma_2 \subset \Sigma_1$, which must bound a compressing disk or a cut disk in $(B_2,\n)$.  But this implies that $\Sigma_1$ is $c$-weakly reducible, so it is stabilized by Lemma \ref{L18}, contradicting the fact that $b_1 = b(K)$.  Thus, $n_1 \geq 4b_1 - 7$, from which it follows that
\[ B(\Sigma) \geq n_1 - b_1 + 1 \geq 4b_1 - 7 - b_1 + 1 = 3b_1 - 6 = 3b(K) - 6.\]
For $b(K) \geq 3$, it holds that $3b(K) -6 \geq b(K) -1$, establishing the lemma.

\end{proof}
\end{lemma}

Suppose $K_{p,q}$ is a torus knot contained in a torus $\Sigma$ such that $S^3 = V \cup_{\Sigma} W$ is a genus one Heegaard splitting.  Let $\A$ and $\n$ be arcs of $K_{p,q}$ such that $K_{p,q}$ is the endpoint union of $\A$ and $\n$.  Leaving $\pd \A = \pd \n$ fixed and pushing $\text{int}(\A)$ slightly into $V$ and $\text{int}(\n)$ slightly into $W$ yields a $(1,1)$-splitting $(V,\A) \cup_{\Sigma} (W,\n)$ of $K_{p,q}$.  We will need another lemma before we calculate the complexity of a torus knot.

\begin{lemma}\label{L20}\cite{morim}
Any torus knot $K_{p,q}$ has a unique $(1,1)$-bridge splitting up to isotopy.
\end{lemma}

Here we classify the pants decompositions that can define one side of a $(1,1)$-splitting.  Suppose $V$ is a solid torus containing one trivial arc $\A$, where $\Sigma = \pd V$ and $\Sigma_{\A} = \Sigma - \pd (\ceta(\A))$.  We claim there are exactly two vertices in $C^*(\Sigma_{\A})$ defining $(V,\A)$.  Let $\Delta$ be a bridge disk for $\A$ in $V$.  There is only one cut disk $C$ for $(V,\A)$, a meridian of $V$ such that $|C \cap \A| = 1$.  There are exactly two compressing disks for $V$, a meridian disk $D_1$ that misses $\A$ and $D_2 = \pd (\ceta(\Delta))$.  Let $u_1 = \pd D_1$, $u_2 = \pd D_2$, and $u_3 = \pd C$.  Since $\pd u_2 \cap \pd u_3 \neq \emp$, $v_1 = \{u_1,u_2\}$ and $v_2 = \{u_1,u_3\}$ are the only two vertices in $C^*(\Sigma_{\A})$ defining $(V,\A)$. \\

Suppose $K_{p,q} = (V,\A) \cup_{\Sigma} (W,\n)$ is the unique $(1,1)$-splitting described above.  Define $u_1' = \pd D_1'$, $u_2' = \pd D_2'$, and $u_3' = \pd C'$ for the disks in $(W,\n)$ corresponding to the definitions for $(V,\A)$ above.  Let $w_1 = \{u_1',u_2'\}$ and $w_2 = \{u_1',u_3'\}$ be the two vertices in $C^*(\Sigma_{\n})$ defining $(W,\n)$.  The arcs $\A,\n$ are isotopic to arcs $\A',\n' \subset \Sigma$, where $\A' \cap \n' = \pd \A' = \pd \n'$, $\A' \cup \n'$ is a $(p,q)$-curve on $\Sigma$, $u_2 = \pd (\Sigma \cap \ceta(\A'))$, and $u_2' = \pd (\Sigma \cap \ceta(\n'))$. \\

We claim that $D(\Sigma) = 3$.  Note that $D(\Sigma) = D(v_i,w_j)$ for some $i,j \in \{1,2\}$.  Since $v_i \cap w_j = \emp$, we immediately have $D(v_i,w_j) \geq 2$.  Suppose by way of contradiction that $D(\Sigma) = 2$.  Then there exists $v_* \in C^*(\Sigma_K)$ such that $v_* \cap v_i \neq \emp$ and $v_* \cap w_j \neq \emp$; hence $v_* = \{u_l,u_k'\}$ for some $l,k \in \{1,2,3\}$.  However, for $l,k \in \{1,3\}$, we have $|u_l \cap u_k'| = 1$.  Additionally, $|u_1 \cap (\A' \cup \n')| = q$, so since $u_1 \cap \A' = \emp$, we have $u_1 \cap \n' \neq \emp$ and thus $u_1 \cap u_2' \neq \emp$.  Similarly, $|u_3 \cap (\A' \cup \n')| = q > 1$ and $|u_3 \cap \A'| = 1$, so $u_3 \cap \n' \neq \emp$ and $u_3 \cap u_2' \neq \emp$.  Similar arguments show that $u_2 \cap u_1'$ and $u_2 \cap u_3'$ are nonempty.  It follows that $D(\Sigma) \geq 3$.  Finally, letting $\gamma \subset \Sigma_K$ denote a $(p,q)$-curve on $\Sigma$ that misses $u_2 \cup u_2'$, we have that $v_1,\{\gamma,u_2\},\{\gamma,u_2'\},w_1$ is a path in $C^*(\Sigma_K)$; thus $D(\Sigma) = 3$ as desired.  As any torus knot $K_{p,q}$ is meridionally small, we have the following:

\begin{theorem}\label{torus}
For any torus knot $K_{p,q}$ with $p > 2$,
\[ B(K_{p,q}) = 2.\]
\begin{proof}
Suppose $(S^3,K_{p,q}) = (V,\A) \cup_{\Sigma} (W,\n)$ is a bridge splitting such that $B(\Sigma) = B(K_{p,q}$), and let $\Sigma_1$ be the critical component of $\Sigma$ with respect to a minimal path between $v$ and $w$ in $C^*(\Sigma_K)$ defining $(V,\A)$ and $(W,\n)$, where $\Sigma_1$ is a $(g_1,b_1)$-bridge surface.  If $g_1 = 0$, then
\[ B(K_{p,q}) = B_0(K_{p,q}) \geq b(K_{p,q}) - 1 \geq 3 - 1= 2\]
by Lemma \ref{L19}.  On the other hand, suppose that $g_1 \geq 2$.  By Lemma \ref{L15},
\[ B(K_{p,q}) \geq 2g_1 + b_1 - 2 \geq 4 + 1 - 2 = 3.\]
Now, suppose that $g_1 = 1$ and $b_1 \geq 2$.  Again, by Lemma \ref{L15},
\[ B(K_{p,q}) \geq 2g_1 + b_1 - 2 \geq 2 + 2 - 2 = 2.\]
Finally, suppose that $g_1 = b_1 = 1$.  By Lemma \ref{L18}, $\Sigma_1$ must be the splitting described above, with $D(\Sigma_1) = 3$ and $B(\Sigma_1) = B(K_{p,q}) = 2$, completing the proof of the theorem.
\end{proof}
\end{theorem}

\begin{figure}[h]
  \centering
    \includegraphics[width=0.9\textwidth]{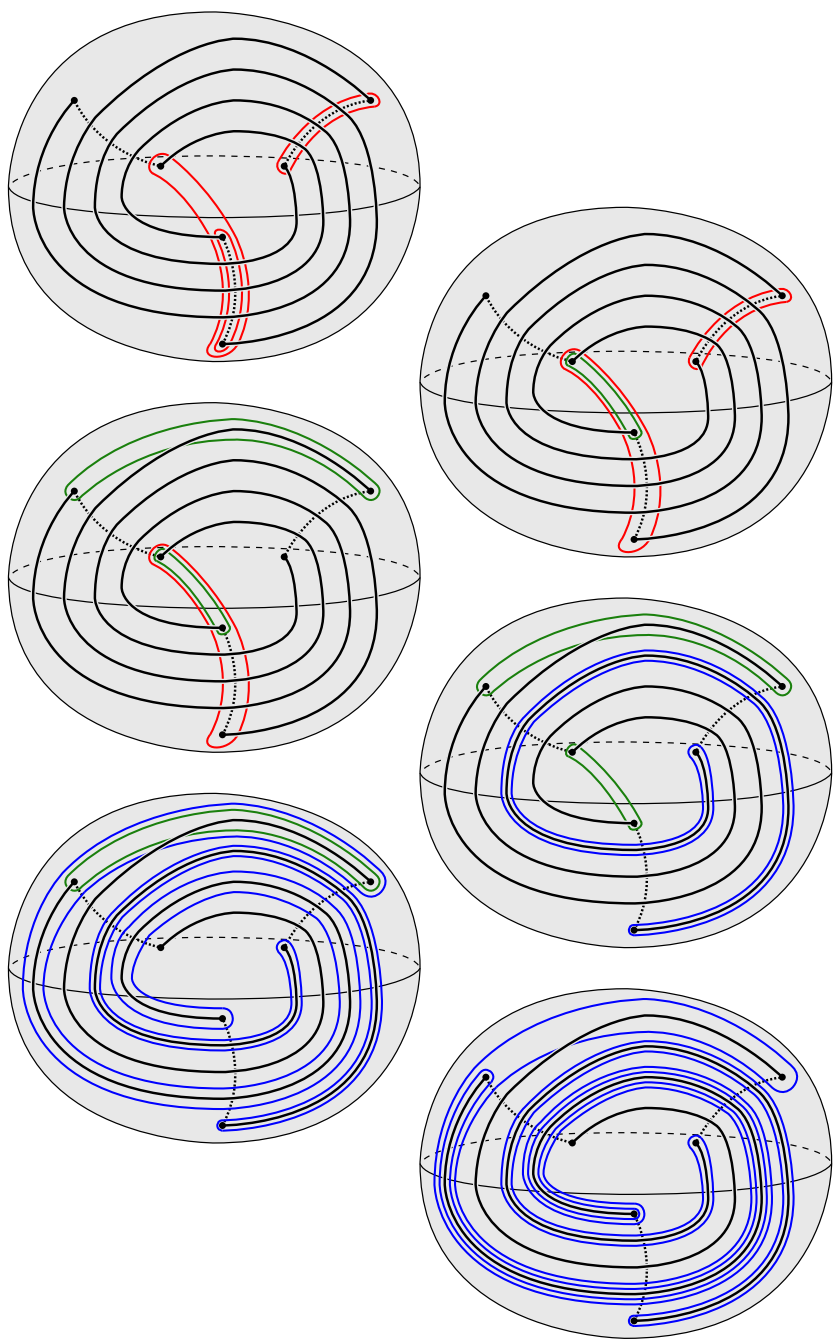}
    \caption{A path in $C^*((\Sigma_0)_K)$, where red curves bound disks in $B_1$ and blue curves bound disks in $B_2$}\label{tor1}
\end{figure}

As an example, consider $K_{3,4}$.  By \cite{schult}, $b(K_{3,4}) = 3$.  Let $(B_1,\A) \cup_{\Sigma_0} (B_2,\n)$ be the $(0,3)$-bridge surface for $K_{3,4}$, which is unique by \cite{ozawa}.  Using the proof of Lemma \ref{L15}, we have that $D(\Sigma_0) \geq 5$.  In Figure \ref{tor1} we exhibit a path of length 5 between vertices $v$ defining $(B_1,\A)$ and $w$ defining $(B_2,\n)$.  Thus, $D(\Sigma_0) = 5$ and $B(\Sigma_0) = 3$. \\

\begin{figure}[h]
  \centering
    \includegraphics[width=0.9\textwidth]{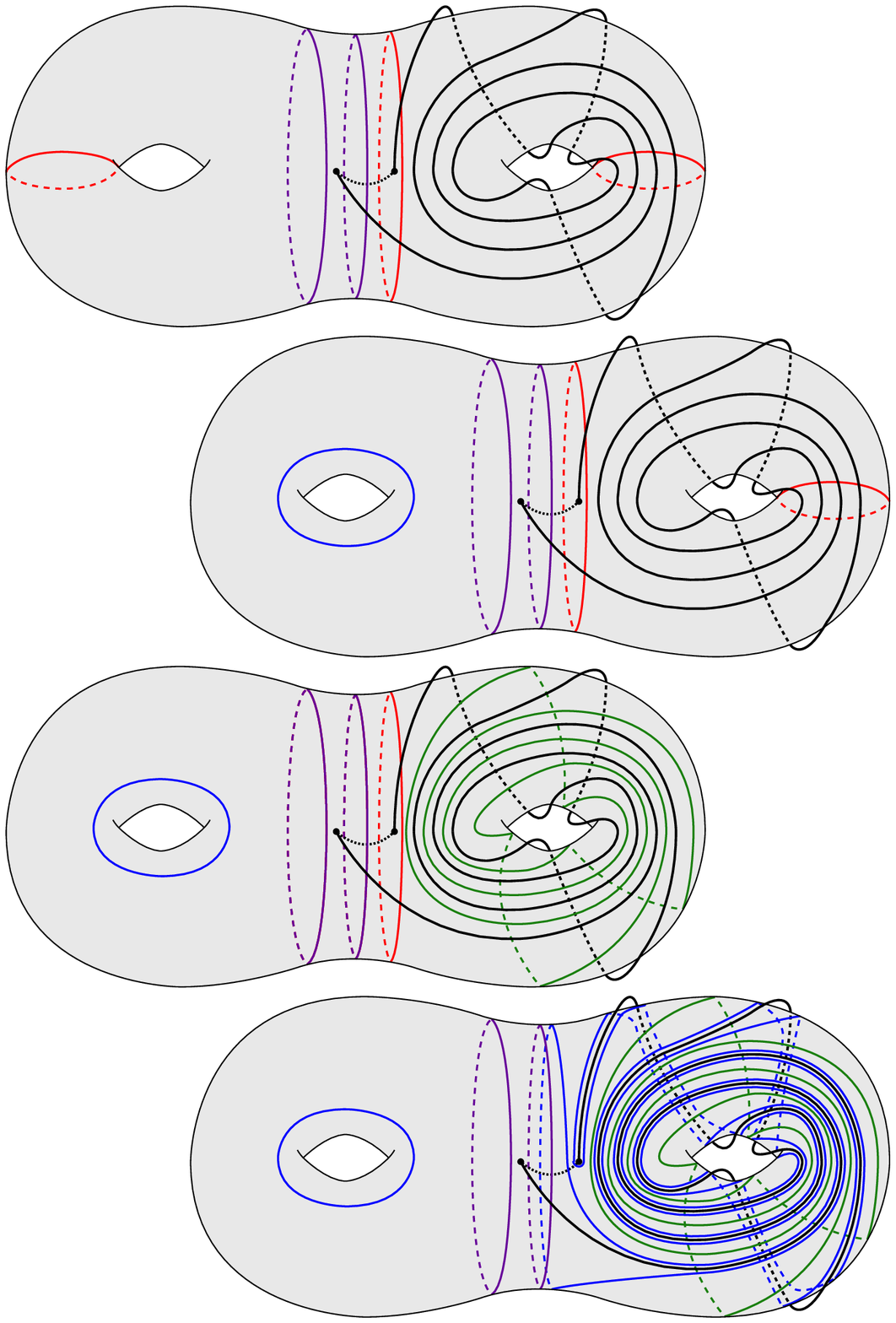}
    \caption{A path (omitting the endpoint $w$) in $C^*((\Sigma_1')_K)$, where red curves bound disks in $V'$, blue curves bound disks in $W$, and purple curves bound disks in both handlebodies}\label{tor2}
\end{figure}

On the other hand, let $\Sigma_1$ be the $(1,1)$-splitting of $K_{3,4}$, and let $(V',\A') \cup_{\Sigma_1'} (W',\n')$ be an elementary stabilization of $\Sigma_1$, so that $\Sigma_1'$ is a $(2,1)$-splitting surface.  By Theorem \ref{torus}, $D(\Sigma_1) = 3$ and $B(\Sigma_1) = 2$, and by Lemma \ref{L7}, $B(\Sigma_1') = B(\Sigma_1) = 2$; hence $D(\Sigma_1') = 4$.  In Figure \ref{tor2}, we exhibit a path of length 4 between vertices $v$ defining $(V',\A')$ and $w$ defining $(W',\n')$ (with a picture of $w$ omitted due to its complexity).

\section{Relationship between Pants Distance and Hyperbolic Volume}\label{hyp}

The literature suggests intriguing evidence of a connection between the pants complex and hyperbolicity.  A compelling piece of this evidence has been produced by Brock, who showed that the pants complex is quasi-isometric to Teichm\"{u}ller space equipped with the Weil-Petersson metric \cite{brock}. \\

In addition, the notion of pants distance of a Heegaard splitting has also been studied by Souto in \cite{souto}.  He uses a slightly different definition than the one that appears in \cite{johnson}: we say for a handlebody $V$ that $v \in P(\pd V)$ \emph{decomposes} $V$ if there exists a collection of curves $v' \subset v$ bounding a collection of compressing disks $\Delta$ in $V$ such that $V - \eta(\Delta)$ is a collection of solid tori.  Then for a Heegaard splitting $M = V \cup_S W$, define
\[ \delta_P(S) = \min \{D^P(v,w): v \text{ decomposes $V$ and $w$ decomposes $W$}\}.\]
Souto proves the following, attributing the result to himself and Brock:
\begin{theorem}\cite{souto}
If $M$ is hyperbolic and $S$ is strongly irreducible, then there exists a constant $L_g$ depending only on $g$ such that
\[ \frac{\delta_P(S)}{L_g} \leq \text{vol}(M) \leq L_g \delta_P(S).\]
\end{theorem}

\begin{figure}[h]
  \centering
    \includegraphics[width=0.55\textwidth]{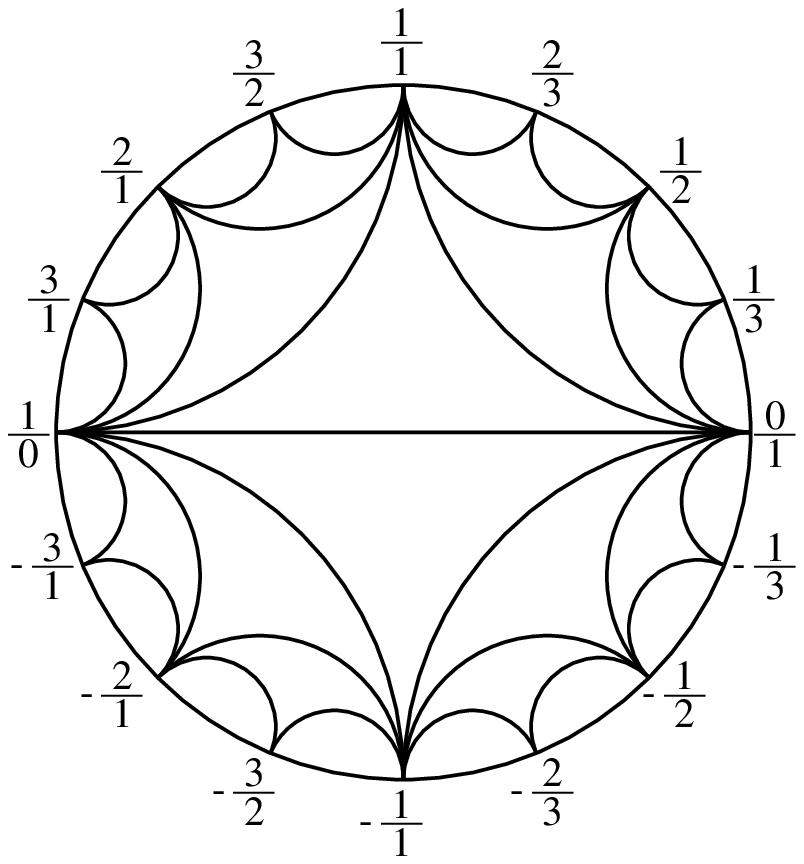}
    \caption{A piece of the Farey graph $\FF$}\label{farey}
\end{figure}

Now, suppose that $K$ is a hyperbolic 2-bridge knot.  We will focus on the pants distance of the $(0,2)$-splitting of $K$, $(S^3,K) = (B_1,\A) \cup_{\Sigma} (B_2,\n)$.  This splitting has surface $\Sigma_K$, a 4-punctured sphere whose pants complex is a Farey graph $\FF$.  The vertex set of $\FF$ is the extended rational numbers $\Q \cup \infty$, where two rational numbers $\frac{a}{b}$ and $\frac{c}{d}$ are joined by an edge if $|ad-bc| = 1$.  The graph $\FF$ can be realized on the hyperbolic plane $\mathbb{H}^2$, with the vertex set contained in the boundary and the edges as geodesics.  For further reference, see \cite{masur}.  A depiction of $\FF$ appears in Figure \ref{farey}. \\

\begin{figure}[h]
  \centering
    \includegraphics[width=.9\textwidth]{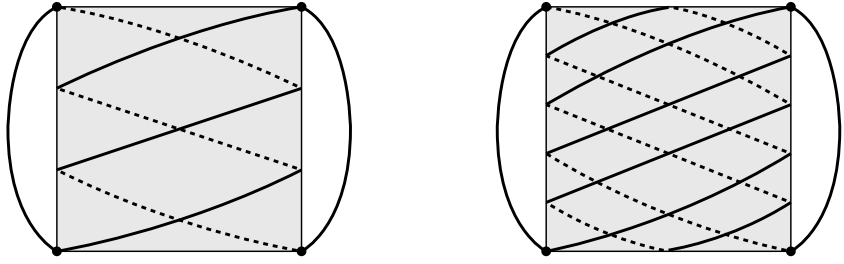}
    \caption{The knots $K_{1/3} = 3_1$ (left) and $K_{2/5} = 4_1$ (right)} \label{ratex}
\end{figure}

The bridge surface $\Sigma_K$ for a 2-bridge knot can be constructed by gluing two unit squares along their boundaries, where the four punctures appear at the vertices.  This is sometimes called a \emph{pillowcase}.  For a reduced rational number $\frac{p}{q}$, we construct a 2-bridge knot by drawing arcs with slope $\frac{p}{q}$ on the surface of the pillowcase and connecting the points $(0,0)$ to $(0,1)$ and $(1,0)$ to $(1,1)$ with trivial arcs outside of the pillowcase.  Pushing the sloped arcs slightly into the ball bounded by the pillowcase yields a $(0,2)$-splitting for the constructed knot, which we call $K_{p/q}$.  See Figure \ref{ratex} for examples. \\

In his classical paper on the subject \cite{schub2}, Schubert proved that two knots $K_{p/q}$ and $K_{p'/q'}$ are isotopic if and only if $q' = q$ and $p' \equiv p \text{ mod } q$ or $pp' \equiv 1 \text{ mod } q$ using the fact that the double branched cover of $S^3$ over $K_{p/q}$ is the lens space $L(q,p)$.  It should also be noted that clearly $K_{-p/q}$ is the mirror image of $K_{p/q}$. \\

Any essential simple closed curve $\gamma$ in $\Sigma_K$ is the boundary of a regular neighborhood of an arc $\A$ connecting distinct punctures, and thus $\gamma$ can be assigned the slope of $\A$, which we denote $m(\gamma)$.  Let $K_{p/q}$ be a 2-bridge knot with $(0,2)$-splitting given by $(S^3,K_{p/q}) = (B_1,\A) \cup_{\Sigma} (B_2,\n)$.  If $v \in P(\Sigma)$ defines $(B_1,\A)$, then $m(v) = \frac{1}{0}$, and if $w \in P(\Sigma)$ defines $(B_2,\n)$, we have $m(w) = \frac{p}{q}$.  Thus, a minimal path between $v$ and $w$ in $P(\Sigma)$ is a minimal path between $\frac{1}{0}$ and $\frac{p}{q}$ in $\FF$.  Any element of $SL(2,\Z)$ induces an automorphism of $\FF$; thus using Schubert's classification and passing to a mirror image if necessary, we may use the automorphism given by
\[ A = \begin{bmatrix} \pm 1 & n \\ 0 & 1 \end{bmatrix} \]
to suppose without loss of generality that $0 < \frac{p}{q} \leq \frac{1}{2}$. \\

A minimal path between $\frac{0}{1}$ and $\frac{p}{q}$ now corresponds uniquely to a continued fraction expansion $[a_1,\dots,a_n]$ of $\frac{p}{q}$, where $a_i > 0$ and $a_1 > 1$.  Successive truncations of this continued fraction yield a path between $\frac{1}{0}$ and $\frac{p}{q}$, although this path is minimal if and only if $a_i > 1$ for all $i$.  If $a_j  =1$ for some $j$, then $[a_1,\dots,a_{j-2}]$ and $[a_1,\dots,a_j]$ are connected by an edge in $\FF$.  For example, $\frac{2}{5} = [2,2]$, and so a minimal path between $\frac{1}{0}$ and $\frac{2}{5}$ is given by
\[ \left\{\frac{1}{0},\frac{0}{1},[2] = \frac{1}{2},[2,2] = \frac{2}{5}.\right\} \]
As a second example, $\frac{3}{11} = [3,1,2]$, and a path between $\frac{1}{0}$ and $\frac{3}{11}$ is given by
\[ \left\{\frac{1}{0},\frac{0}{1},[3] = \frac{1}{3},[3,1] = \frac{1}{4},[3,1,2] = \frac{3}{11},\right\}\]
but a minimal path omits $\frac{1}{2}$.  It follows immediately that the length of a minimal path is at least one more than half the length of the corresponding continued fraction expansion; hence if $\frac{p}{q} = [a_1,\dots,a_n]$, then
\begin{equation}\label{pineq}
\frac{n}{2} \leq D^P(\Sigma) - 1 \leq n.
\end{equation}
For a more detailed discussion of 2-bridge knots and paths in the Farey graph, refer to \cite{hatchthurst}. \\

\begin{figure}[h]
  \centering
    \includegraphics[width=.9\textwidth]{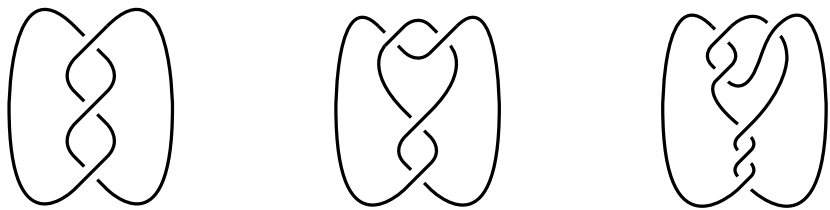}
    \caption{The knots $[3] = 3_1$ (left), $[2,2] = 4_1$ (middle), and $[3,1,2] = K_{3/11}$ (right)} \label{tang}
\end{figure}

The continued fraction expansion $\frac{p}{q} = [a_1,\dots,a_n]$ also corresponds to a rational tangle diagram $D_{p,q}$ of $K_{p,q}$, where each $a_i$ corresponds to some number of horizontal or vertical twists.  The assumption $a_i > 0$ ensures that $D_{p,q}$ is alternating.  We will not go into detail in this article, but a good reference for the interested reader is \cite{kauflamp}.  Note that \cite{kauflamp} uses numerator closure of a rational tangle to construct $K_{p,q}$, whereas we follow the convention of \cite{hatchthurst} in using the denominator closure of the tangle.  See Figure \ref{tang} for examples. \\

A \emph{twist region} in a knot diagram $D$ is a connected sequence of bigons in $\R^2 \setminus D$, where two bigons are adjacent if they share a vertex.  Two crossings are said to be \emph{twist equivalent} if they appear as vertices in the same twist region, and the number of equivalence classes is the \emph{twist number} $tw(D)$ of the diagram.  As $[a_1,\dots,a_n,1] = [a_1,\dots,a_n+1]$, we may always pass to a continued fraction $[a_1,\dots,a_n]$ with $a_n > 1$.  Thus, the diagram $D_{p,q}$ corresponding to $[a_1,\dots,a_n]$ satisfies $tw(D_{p,q}) = n$. \\

The following theorem of Lackenby reveals a remarkable connection between the twist number of certain diagrams of a knot $K$ and the hyperbolic volume of $E(K)$ (see \cite{lackenby} for the definitions of \emph{twist-reduced} and \emph{prime}):

\begin{theorem} \cite{lackenby}
Suppose $D$ is a twist-reduced prime alternating diagram of a hyperbolic link $K$ and $v_3$ is the volume of a regular hyperbolic ideal 3-simplex.  Then
\[ v_3 (tw(D) - 2) \leq \text{vol}(E(K)) < 16v_3(tw(D) - 1).\]
\end{theorem}
In an appendix to this paper, Agol and D. Thurston improve the upper bound, showing that in fact
\begin{equation}\label{pineq2}
v_3 (tw(D) - 2) \leq \text{vol}(E(K)) < 10v_3(tw(D) - 1)
\end{equation}
and exhibiting a family of links $K_i$ with diagrams $D_i$ whose volumes approach $10v_3(tw(D_i) -1)$ asymptotically as $i \rightarrow \infty$.  We note that the upper bound is true for all diagrams. \\

Using the fact the diagram $D_{p/q}$ is prime, alternating, and twist-reduced, we may combine the inequalities (\ref{pineq}) and (\ref{pineq2}) to arrive at the following theorem:

\begin{theorem}\label{bounds}
Suppose $K$ is a hyperbolic 2-bridge knot, $\Sigma$ is a $(0,2)$-splitting surface for $K$, and $v_3$ is the volume of a regular hyperbolic ideal 3-simplex.  Then
\[ v_3(D^P(\Sigma) - 3) \leq \text{vol}(E(K)) < 10v_3(2D^P(\Sigma) - 3).\]
\end{theorem}

\section{Questions}\label{quest}

It follows as a consequence of Lemma \ref{L15} and the discussion in Sections \ref{calc} and \ref{hyp} that $B^P(K) = 1$ if and only if $K$ is the torus knot $K_{2,n}$ for some $n$.  However, this and the other low complexity calculations of Theorem \ref{torus} are the only complexities we have computed.

\begin{question}
Are there reasonable computations of $B(M,K)$ or $B^P(M,K)$ for certain families of knots and manifolds?
\end{question}

One of the natural properties of Hempel's distance of a Heegaard splitting is that as a consequence of Haken's Lemma, the distance of a splitting surface $\Sigma$ of a 3-manifold $M$ is zero whenever $M$ contains an essential sphere or $\Sigma$ is stabilized.  Thus, we would hope that for a bridge splitting surface $\Sigma_K$, the distances $D(\Sigma_K)$ and $D^P(\Sigma_K)$ recognize the topology of $M$, $K$, and $\Sigma_K$.  More specifically, we ask the following question:

\begin{question}
Suppose that $(M,K) = (V,\A) \cup_{\Sigma} (W,\n)$ is a bridge splitting, and $\gamma$ is a simple closed curve that bounds compressing or cut disks in both $(V,\A)$ and $(W,\n)$.  Let $v = v_0,\dots,v_n = w$ be a minimal path in $C^*(\Sigma_K)$ or $P(\Sigma_K)$ between vertices $v$ and $w$ defining $(V,\A)$ and $(W,\n)$.  Is $\gamma \in v \cap w$?  If so, is $\gamma \in v_i$ for all $i$?
\end{question}
An answer to this question may have implications about the behavior of either complexity under connected sums. \\

In the example at the end of Section \ref{calc}, we compute several splitting complexities for the knot $K_{3,4}$.  Note that for the $(0,3)$-splitting surface $\Sigma_0$ of $K_{3,4}$, we have $B(\Sigma_0) > B(K)$, while for the $(1,1)$-splitting surface $\Sigma_1$, $B(\Sigma_1) = B(K)$.  Thus $\Sigma_1$ is the critical component of infinitely many splittings $\Sigma$ such that $B(\Sigma) = B(K)$, but $\Sigma_0$ is never the critical component of  such a splitting.  In the case of torus knots, the critical component of infinitely many splittings is a twice-punctured torus, which inspires the following question:

\begin{question}
Suppose that the $(g_1,b_1)$-splitting surface $\Sigma_1$ for a knot $K$ is the critical component of some splitting $\Sigma$.  Does the pair $(g_1,b_1)$ carry any information about the topology of $K$?  In particular, under what conditions is $g_1$ necessarily the minimal genus Heegaard surface $S$ such that $K \subset S$?
\end{question}
The genus of such $S$ is called the $h$-genus of $K$; see \cite{morim2}. \\

The strength of the approach of Theorem \ref{bounds} is that it does not depend on a specific knot diagram; instead it uses information about compressing and cut disks on either side of a bridge splitting to obtain information about volume.  An extension of Theorem \ref{bounds} would be of interest:

\begin{question}
Are there linear functions of $D^P(\Sigma)$ which yield upper and lower bounds for $\text{vol}(E(K))$ for knots $K$ that are not 2-bridge?  What about $K$ non-alternating?
\end{question}

\begin{question}
Is there a linear function of $D^P(\Sigma)$ which yields lower bounds for $\text{vol}(E(K))$ that are asymptotically sharp for some family of knots?
\end{question}

\end{document}